\documentclass[a4paper,11pt]{article}

\usepackage{mathrsfs}
\usepackage{amsmath,amscd}
\usepackage{amssymb, amsmath}
\usepackage{graphicx}
\usepackage{color}

\allowdisplaybreaks

\newcommand{\newcom}{\newcommand}

\newcom{\cA}{{\mathcal A}}
\newcom{\cB}{{\mathcal B}}
\newcom{\cC}{{\mathcal C}}
\newcom{\cD}{{\mathcal D}}
\newcom{\cE}{{\mathcal E}}
\newcom{\cF}{{\mathcal F}}

\newcom{\cJ}{{\mathcal J}}

\newcom{\cL}{{\mathcal L}}
\newcom{\cM}{{\mathcal M}}
\newcom{\cP}{{\mathcal P}}
\newcom{\cS}{{\mathcal S}}
\newcom{\cQ}{{\mathcal Q}}
\newcom{\cT}{{\mathcal T}}
\newcom{\cY}{{\mathcal Y}}
\newcom{\cZ}{{\mathcal Z}}
\newcom{\R}{\mathbb R}
\newcom{\T}{\mathbb T}
\newcom{\N}{\mathbb N}
\newcom{\Z}{\mathbb Z}
\newcom{\C}{\mathbb C}
\newcom{\E}{\mathbb E}

\newcom{\e}{\epsilon}

\newcom{\al}{\alpha}
\newcom{\be}{\beta}
\newcom{\del}{\delta}

\newcom{\ga}{\gamma}
\newcom{\Ga}{\Gamma}

\newcom{\Lam}{\Lambda}
\newcom{\lam}{\lambda}

\newcom{\Om}{\Omega}
\newcom{\om}{\omega}

\newcom{\Si}{\Sigma}
\newcom{\si}{\sigma}
\newcom{\s}{\varsigma}

\newcom{\tht}{\theta}
\newcom{\dtri}{\nabla}
\newcom{\tri}{\triangle}

\newcom{\f}{\frac}
\newcom\na{\nabla}
\newcom{\Del}{\Delta}
\newcom{\ep}{\epsilon}

\newcom{\p}{\partial}

\newcom{\uep}{{\bf u}_{\epsilon}}
\newcom{\nep}{n_{\epsilon}}
\newcom{\cep}{c_{\epsilon}}

\newcom{\beq}{\begin{equation}}
\newcom{\eeq}{\end{equation}}
\newcom{\ben}{\begin{eqnarray}}
\newcom{\een}{\end{eqnarray}}
\newcom{\beno}{\begin{eqnarray*}}
\newcom{\eeno}{\end{eqnarray*}}
\newcom{\bal}{\begin{aligned}}
\newcom{\eal}{\end{aligned}}

\title{Stability of Big Solitons in a Competitive Power Nonlinear
Schr\"odinger Equation}

\author{Jian Zhang ,
\quad Mengxue Bai \thanks{Corresponding author. E-mail: mengxuebai@163.com}\\
\emph{\small  School of Mathematical Sciences, University of Electronic Science and Technology of China, }\\
\emph{\small Chengdu 611731, China}\\
  }

\date{}

\begin{document}
\maketitle
{\small
{\textbf{Abstract:} By introducing and solving two correlative constrained variational problems as well as spectrum analysis, an approach to fix soliton  frequency from the
prescribed mass for nonlinear Schr\"odinger equations is found,  and  an open problem in normalized solutions is answered. Then existence and orbital stability of big solitons depending on frequencies for nonlinear Schr\"odinger equation with competitive power nonlinearity is proved for the first time.  In addition multi-solitons of the equation with different speeds are constructed by stable big solitons.
\\

{\textbf{Mathematics Subject Classification (2010):} 35Q55; 35J50; 37K40.

{\textbf{Keywords:} {nonlinear Schr\"odinger equation; variational method; normalized solution; spectrum analysis;   stability of soliton}

\section{Introduction}
\renewcommand{\theequation}{\thesection.\arabic{equation}}
\newtheorem{Definition}{Definition}[section]
\newtheorem{Theorem}{Theorem}[section]
\newtheorem{Proposition}{Proposition}[section]
\newtheorem{Lemma}{Lemma}[section]
\newtheorem{Remark}{Remark}[section]
\newtheorem{Corollary}{Corollary}[section]
\numberwithin{equation}{section}

\par   Consider the competitive power nonlinear Schr\"odinger equation
\begin{align}
i\varphi_{t}+\Delta \varphi+\rvert \varphi\lvert^{\frac{4}{d}}\varphi-\rvert\varphi\lvert^{p-1}\varphi=0,\;\;\;\;\;(t,x)\in\ \mathbb{R} \times \mathbb{R}^{d}
\end{align}
with $d\geq 2$ and $1+\frac{4}{d}<p<\frac{d+2}{(d-2)^{+}}$. Here and hereafter we denote $\frac{d+2}{(d-2)^{+}}=\infty$ for $d=2$ and $\frac{d+2}{(d-2)^{+}}=\frac{d+2}{d-2}$ for $d\geq 3$. We also denote $\int_{\mathbb{R}^{d}}\cdot dx$ by $\int\cdot dx$, $L^{2}(\mathbb{R}^{d})$ by $L^{2}$ and $H^{1}(\mathbb{R}^{d})$ by $H^{1}$.
\par  We are interested in  solitons $e^{i\omega t}u(x)$ of (1.1), where $\omega\in\mathbb{R}$ is called frequency of the soliton and $u(x)$ satisfies the time-independent nonlinear Schr\"odinger equation
\begin{align}
\Delta u+\rvert u\lvert^{\frac{4}{d}}u-\rvert u\lvert^{p-1}u-\omega u=0,\;\;\;\;\;u\in\ H^{1}.
\end{align}
\par  (1.1) is proposed and studied by many papers, such as Tao, Visan and Zhang \cite{TVZ2007} for global well-posedness and scattering; Soave \cite{S2020} for normalized solutions; Le Coz, Martel and Rapha\"el \cite{CMR2016} for minimal mass blow up solutions; Fukaya and Hayashi \cite{FH2021} for instability of solitons et al. Classic results of (1.2) come from Berestycki and Lions \cite{BL1983}; Gidas, Ni and Nirenberg \cite{GNN1979}; McLeod \cite{M1993}; Strauss \cite{S1977}. A solution of (1.2) satisfying prescribed mass
\begin{align}
\int\rvert u\lvert^{2}dx=m,\;\;\;\;m > 0
\end{align}
is called normalized solution. Besides motivations in mathematical physics, normalized solutions are also of interest in the framework of ergodic Mean Field Games systems \cite{CV2017}. 
Recent studies on normalized solutions refer to \cite{BJS2016, BMRV2021, BZZ2020, S2020} et al. A crucial problem of normalized solutions is whether the mapping from prescribed mass $m$ to soliton frequency $\omega$ is injective, which still remains extensively open except for the equation with single power nonlinearity (see \cite{CL1982}). And this problem is also crucial for stability of solitons depending on frequencies instead of mass (see \cite{GSS1987, GSS1990, S2020, Z2000}).
\par  According to Berestycki and Lions \cite{BL1983} as well as McLeod \cite{M1993}, Fukuizumi \cite{F2003} and Kawano \cite{K2011} got that (1.2) possesses a unique positive radial solution $Q_{\omega}(x)$ if and only if $\omega\in (0,\;\omega_{p})$, where
\begin{align}
\omega_{p}=\frac{2(p-1-4/d)}{(2+4/d)(p-1)}[\frac{(p+1)4/d}{(2+4/d)(p-1)}]^{(4/d)/(p-1-4/d)}.
\end{align}
\par  From Weinstein \cite{W1983} and Kwong \cite{K1989}, the nonlinear scalar field equation
\begin{align}
\Delta u+\rvert u\lvert^{\frac{4}{d}}u-u=0,\;\;\;\;\;u\in\ H^{1}
\end{align}
possesses a unique positive radial solution denoted by $q$. 
The Hamilton functional related to (1.1) is that
\begin{align}
E(u)=\int\rvert \nabla u\lvert^{2}dx-\frac{1}{1+\frac{2}{d}}\int\rvert u\lvert^{2+\frac{4}{d}}dx+\frac{2}{p+1}\int\rvert u\lvert^{p+1}dx,
\end{align}
where $u\in H^{1}$. For $u\in H^{1}$, the Pohozaev functional related to (1.2) is that
\begin{align}
I(u)=2\int\rvert\nabla u\lvert^{2}dx-\frac{2d}{d+2}\int\rvert u\lvert^{2+\frac{4}{d}}dx+\frac{p-1}{p+1}d\int\rvert u\lvert^{p+1}dx.
\end{align}
\par  For $\int q^{2}dx<m<\infty$, we set the Cazenave-Lions type constrained variational problem (see \cite{CL1982})
\begin{align}
d_{m}= inf_{\{u\in H^{1},\;\int \rvert u\lvert^{2}dx=m\}}E(u).
\end{align}
We will prove that (1.8) possesses a minimizer $Q_{m}\in H^{1}$. Moreover there exists a Lagrange multiplier $\omega_{m}\in\mathbb{R}$ such that $Q_{m}$ is the unique positive radial solution of (1.2) with $\omega=\omega_{m}$. Therefore we can define
\begin{align}
\omega_{q}=inf\{\omega\in\mathbb{R}\big\rvert\;\int Q^{2}_{\omega}dx > \int q^{2}dx\}.
\end{align}
and show that $0 \leq \omega_{q} < \omega_{p}$, where $Q_{\omega}$ is the positive radial solution of (1.2). Now for $\omega\in (\omega_{q}, \omega_{p})$ and $\int q^{2}dx<m<\infty$, we introduce the following constrained variational problem
\begin{align}
d_{\omega}= inf_{\{u\in H^{1},\;\int q^{2}dx < \int \rvert u\lvert^{2}dx \leq m,\;I(u)=0\}}[E(u)+\omega \int \rvert u\lvert^{2}dx].
\end{align}
(1.10) is a multi-constrained variational problem. We can prove that (1.10) is solvable. (1.10) with (1.8) provides a new correlative variational framework, by which we can fix frequency from the prescribed mass for nonlinear Schr\"odinger equations. Combining with spectral analysis, we claim the following theorem.

\noindent\textbf{Theorem A. } A one-to-one mapping between $m\in (\int q^{2}dx,\;\infty)$ and $\omega\in (\omega_{q},\;\omega_{p})$ is determined by $m=\int Q^{2}_{\omega}dx$, where $Q_{\omega}$ is the unique positive solution of (1.2) with $\omega\in (\omega_{q},\;\omega_{p})$. Moreover for $\omega\in (\omega_{q},\;\omega_{p})$,
\begin{align}
\frac{dm}{d\omega}=\frac{d}{d\omega}\int Q_{\omega}^{2}dx > 0.
\end{align}

\par  We see that Theorem A establishes a one-to-one mapping from prescribed mass to frequency in normalized solutions. The approach introduced in the present paper can be used to deal with more nonlinear Schr\"odinger equations, such as general double power nonlinear Schr\"odinger equations (see \cite{CMR2016} and \cite{S2020}), nonlinear Hartree equations (see \cite{CL1982}), nonlinear Schr\"odinger equations with potentials (see \cite{BMRV2021}, \cite{Z2000} and \cite{Z2005} ) et al. We provide an approach to answer the aforementioned open problem in normalized solutions.
\par  In the proof of Theorem A, we synchronously get stability of solitons depending on frequencies.
\\
\noindent\textbf{Theorem B. } For $\omega\in (\omega_{q},\;\omega_{p})$, there exists a big soliton $e^{i\omega t}Q_{\omega}(x)$ of (1.1) satisfying $\int Q^{2}_{\omega}(x)dx > \int q^{2}dx$, where $Q_{\omega}$ is the unique positive solution of (1.2). Moreover the big soliton $e^{i\omega t}Q_{\omega}(x)$ with $\omega\in (\omega_{q},\;\omega_{p})$ is orbitally stable.
\par  Since the stable soliton $e^{i\omega t}Q_{\omega}(x)$ in Theorem B satisfies $\int Q^{2}_{\omega}(x)dx > \int q^{2}dx$, we call it big soliton of (1.1). We see that Theorem A may directly deduce Theorem B (\cite{GSS1987}). But indeed the proof of Theorem A in the present paper depends on Theorem B.
Stability of solitons is a crucial topic in understanding the dynamics of nonlinear dispersive evolution equations \cite{T2009}. We recall general nonlinear Schr\"odinger equation
\begin{align}
i\varphi_{t}+\Delta \varphi+f(\rvert\varphi\lvert^{2})\varphi=0.
\end{align}
When $f(\rvert\varphi\lvert^{2})\varphi=\rvert\varphi\lvert^{p-1}\varphi$ with $1+\frac{4}{d}\leq p < \frac{d+2}{(d-2)^{+}}$, that is mass critical or supercritical nonlinear case, the soliton $e^{i\omega t}u(x)$ with $\omega > 0$ is unstable to blow up (see \cite{BC1981}, \cite{W1983}). When
$f(\rvert\varphi\lvert^{2})\varphi=\rvert\varphi\lvert^{p-1}\varphi$ with $1 < p < 1+\frac{4}{d}$, that is mass subcritical nonlinear case, the soliton $e^{i\omega t}u(x)$ with $\omega > 0$ is orbitally stable (see \cite{CL1982}). When (1.12) has no scaling invariance, it becomes very difficult to prove stability of solitons depending on frequencies (\cite{GSS1987}, \cite{GSS1990}). In this case, for $f(\rvert\varphi\lvert^{2})\varphi$ subject to mass subcritical nonlinearity, Weinstein \cite{W1986} proved stability of solitons depending on frequencies by spectrum analysis. In Theorem B, for $f(\rvert\varphi\lvert^{2})\varphi$ subject to mass critical and supercritical nonlinearity, we get existence and stability of big solitons depending on frequencies by the new correlative variational framework introduced in the present paper. 
 It is known that stability  is   substantially sought in global dynamics of nonlinear Schr\"odinger equations (see \cite{NS2012, S2009}).
\par  From Theorem B, in terms of the bootstrap argument and compactness method introduced by Martel, Merle and Tsai (see \cite{MM2006}, \cite{MMT2006} and \cite{M1990}), we construct multi-solitons of (1.1) by stable big solitons.
\par \noindent\textbf{Theorem C. } For $K \geq 2$ and $k=1,\;2,\;\cdot\cdot\cdot,\;K$, taking $\omega_{k}\in (\omega_{q},\;\omega_{p})$, $\gamma_{k}\in \mathbb{R}$, $x_{k}\in \mathbb{R}^{d}$, $v_{k}\in \mathbb{R}^{d}$ with $v_{k}\neq v_{k'}$ to $k\neq k'$ and
\begin{eqnarray}
R_{k}(t,x)=Q_{\omega_{k}}(x-x_{k}-v_{k}t)e^{i(\frac{1}{2}v_{k}x-\frac{1}{4}\rvert v_{k}\lvert^{2}t+\omega_{k}t+\gamma_{k})},
\end{eqnarray}
there exists a solution $\varphi(t,\;x)$ of (1.1) such that
\begin{align}
\lim\limits_{t \rightarrow +\infty}\|\varphi(t)-\sum_{k=1}^{K}R_{k}(t)\|_{H^{1}}=0.
\end{align}
The solution $\varphi(t,\;x)$ of (1.1) holding (1.14) is called multi-soliton of (1.1).
\par  Multi-solitons are concerned with the famous soliton resolution conjecture \cite{T2009}. According to stability of solitons obtained in \cite{CL1982} and \cite{W1986}, Martel and Merle constructed multi-solitons with different speeds for (1.12) to mass subcritical nonlinearity  and $f(\rvert\varphi\lvert^{2})\varphi$ subject to mass subcritical nonlinearity. Merle \cite{M1990}, C$\hat{o}$te, Martel and Merle \cite{CMM2011} constructed multi-solitons with different speeds by unstable solitons for (1.12) subject to mass critical and supercritical nonlinearities. C$\hat{o}$te and Le Coz \cite{CC2011} constructed high speed excited multi-solitons of (1.12) by unstable solitons. In Theorem C, we get  multi-solitons with different speeds constructed  by stable big solitons for $f(\rvert\varphi\lvert^{2})\varphi$ subject to mass critical and supercritical nonlinearities. Soliton resolution conjecture remains still open for nonlinear Schr\"odinger equations, which depends on stability of multi-solitons of (1.12) (see \cite{MMT2006}). In general, multi-solitons constructed by unstable solitons are unstable (see \cite{CC2011}). Thus, the existence of multi-solitons constructed by stable solitons is  the first  step to realize soliton resolution conjecture.

\par This paper is organized as follows. In section 2, we show global well-posedness of the Cauchy problem for (1.1) in the energy space as well as the sufficient and necessary conditions for (1.2) to possess a unique positive radial solution. In section 3, we introduce and solve two correlative constrained variational problems. In section 4, we prove orbital stability of big solitons of (1.1) depending on frequencies and complete the proofs of Theorem A and Theorem B. In section 5,  we construct multi-solitons with different speeds of (1.1) by stable big solitons.

\section{Well-Posedness}
For $t_{0}\in \mathbb{R}$, we impose the initial data of (1.1) as follows.
\begin{align}
\varphi(t_{0},x)=\varphi_{0}(x),\;\;\;x\in \mathbb{R}^{d},\;\;\;d\geq 2.
\end{align}
In $H^{1}$, we define the energy functional
\begin{align}
E(\varphi):=\int\rvert \nabla \varphi\lvert^{2}dx-\frac{1}{1+\frac{2}{d}}\int\rvert\varphi\lvert^{2+\frac{4}{d}}dx+\frac{2}{p+1}\int\rvert\varphi\lvert^{p+1}dx;
\end{align}
the mass functional
\begin{align}
M(\varphi):=\int\rvert\varphi\lvert^{2}dx;
\end{align}
the momentum functional
\begin{align}
P(\varphi):=Im\int \overline{\varphi}\nabla \varphi dx.
\end{align}
\par \noindent\textbf{Theorem 2.1.} Let $d\geq 2$, $1+\frac{4}{d} < p  < \frac{d+2}{(d-2)^{+}}$ and $\varphi_{0}\in H^{1}$. Then the Cauchy problem (1.1)-(2.1)
possesses a unique global solution $\varphi(t,x)\in C(\mathbb{R},\;H^{1})$ with mass conservation $M(\varphi)=M(\varphi_{0})$, energy conservation $E(\varphi)=E(\varphi_{0})$ and momentum conservation $P(\varphi)=P(\varphi_{0})$ for all $t\in \mathbb{R}$.
\\
\noindent\textbf{Proof.} By \cite{C2003}, for $\varphi_{0} \in H^{1}$, there exists a unique solution $\varphi (t,\;x)$ of the Cauchy problem (1.1)-(2.1) in $C((-T,\;T);\;H^{1})$ to some $T>0$ (maximal existence time). And $\varphi(t,\;\cdot)$ satisfies mass conservation ~~$M(\varphi)=M(\varphi_{0})$, energy conservation ~~$E(\varphi)=E(\varphi_{0})$ and momentum conservation ~~$P(\varphi)=P(\varphi_{0})$ for
all $t\in (-T,\;T)$. Furthermore one has the alternatives: $T=\infty$ (global existence) or else $T<\infty$ and $lim_{t\rightarrow T}\|\nabla\varphi\|_{L^{2}}=\infty$ (blow up).
\par Since $1+\frac{4}{d} < p  < \frac{d+2}{(d-2)^{+}}$, from the interpolation inequality for $\varphi\in H^{1}$, there exists $0 < \theta < 1$ such that $\frac{1}{2+\frac{4}{d}}=\frac{\theta}{p+1}+\frac{1-\theta}{2}$ and
\begin{align}
\|\varphi\|_{L^{2+\frac{4}{d}}} \leq \|\varphi\|^{1-\theta}_{L^{2}}\cdot \|\varphi\|^{\theta}_{L^{p+1}}.
\end{align}
By Young inequality, it follows that for arbitrary $\varepsilon > 0$ there exists a number $C(\varepsilon,\;p,\;\|\varphi\|_{L^{2}}) > 0$ depending on $\varepsilon,\;p$ and $\|\varphi\|_{L^{2}}$ such that
\begin{align}
\frac{1}{1+\frac{2}{d}}\int\rvert\varphi\lvert^{2+\frac{4}{d}}dx\leq C(\varepsilon,\;p,\;\|\varphi\|_{L^{2}})+\varepsilon \int\rvert\varphi\lvert^{p+1}dx.
\end{align}
Take $0 < \varepsilon < \frac{2}{p+1}$, we have that there exists a number $C(\varepsilon,\;d,\;p,\;\|\varphi\|_{L^{2}}) > 0$ depending on $\varepsilon,\;d,\;p$ and $\|\varphi\|_{L^{2}}$, such that
\begin{align}
E(\varphi)&=\int\rvert\nabla \varphi\lvert^{2}dx-\frac{1}{1+\frac{2}{d}}\int\rvert\varphi\lvert^{2+\frac{4}{d}}dx+\frac{2}{p+1}\int\rvert\varphi\lvert^{p+1}dx\notag\\
&\geq \int\rvert\nabla \varphi\lvert^{2}dx-C(\varepsilon,\;d,\;p,\;\|\varphi\|_{L^{2}})+(\frac{2}{p+1}-\varepsilon)\int\rvert\varphi\lvert^{p+1}dx.
\end{align}
From mass and energy conservations, it follows that
\begin{align}
\int\rvert\nabla \varphi\lvert^{2}dx \leq E(\varphi_{0})+C(\varepsilon,\;d,\;p,\;\|\varphi_{0}\|_{L^{2}}),
\end{align}
where $C(\varepsilon,\;d,\;p,\;\|\varphi_{0}\|_{L^{2}})$ is a positive constant number depending on $\varepsilon,\;d,\;p$ and $\|\varphi_{0}\|_{L^{2}}$. It yields that $\|\nabla\varphi\|_{L^{2}}$ is bounded for $t\in (-T,\;T)$ with any $T<\infty$. Therefore we get that $\varphi(t,\;x)$ globally exists in $t\in (-\infty,\;\infty)$. Moreover, mass conservation, energy conservation and momentum conservation are true to all $t\in \mathbb{R}$.
\par  This proves Theorem 2.1.
\par  By Weinstein \cite{W1983} and Kwong \cite{K1989}, we state the following lemma.
\par \noindent\textbf{Lemma 2.2.}  (1.5) possesses a unique positive radial solution $q=q(x)$. Moreover for $d \geq 2$ and $\phi\in H^{1}$, one has that
\begin{align}
\int\rvert\phi\lvert^{2+\frac{4}{d}}dx\leq \frac{2+d}{d}(\int q^{2}dx)^{-\frac{2}{d}}(\int\rvert\nabla \phi\lvert^{2}dx)(\int\rvert\phi\lvert^{2}dx)^{\frac{2}{d}}.
\end{align}
\par \noindent\textbf{Theorem 2.3.} Let $d\geq 2$, $1+\frac{4}{d} < p < \frac{d+2}{(d-2)^{+}}$. Then the necessary conditions for the nonlinear elliptic equation
\begin{align}
\Delta u-\omega u+u\rvert u\lvert^{\frac{4}{d}}-u\rvert u\lvert^{p-1}=0,\;\;u\in H^{1}
\end{align}
to possess non-trivial solutions are $\omega > 0$ and $\int \rvert u\lvert^{2}dx > \int q^{2}dx$, where $q$ is the unique positive radial solution of (1.5) and $u$ is the non-trivial solution of (2.10).
\par \noindent\textbf{Proof.} Let $u(x)$ be a non-trivial solution of (2.10). From the Pohozaev's identity (see \cite{P1965}), we have that
\begin{align}
\int(\frac{2d}{p+1}+2-d)\rvert \nabla u\lvert^{2}+(\frac{d}{1+\frac{2}{d}}-\frac{2d}{p+1})\rvert u\lvert^{2+\frac{4}{d}}+(\frac{2d}{p+1}-d)\omega\rvert u\lvert^{2}dx=0,
\end{align}
\begin{align}
2\int \rvert\nabla u\lvert^{2}dx-\frac{2d}{d+2}\int\rvert u\lvert^{2+\frac{4}{d}}dx+\frac{(p-1)d}{p+1}\int\rvert u\lvert^{p+1}dx=0.
\end{align}
Since $d\geq 2$, $1+\frac{4}{d} < p < \frac{d+2}{(d-2)^{+}}$, from (2.11) it follows $\omega> 0$.
By (2.9) and (2.12), it follows that $\int \rvert u\lvert^{2}dx > \int q^{2}dx$.
\par  This proves Theorem 2.3.
\par  According to Gidas, Ni and Nirenberg \cite{GNN1979}, every positive solution of (2.10) is radially symmetric. From Berestycki and Lions \cite{BL1983}, a positive radial solutions of (2.10) exists if and only if $\omega > 0$ and there exists $\zeta > 0$
such that
\begin{align}
\int^{\zeta}_{0}(\rvert s\lvert^{\frac{4}{d}}s-\rvert s\lvert^{p-1}s-\omega s)ds > 0.
\end{align}
Then by Fukuizumi \cite{F2003}, Kawano \cite{K2011} claims the following lemma.
\par \noindent\textbf{Lemma 2.4.} Let $d \geq 2$ and $1+\frac{4}{d} < p < \frac{d+2}{(d-2)^{+}}$. Then (2.10) possesses a unique positive radial solution $Q_{\omega}(x)$ if and only if $\omega\in (0,\;\omega_{p})$, where $\omega_{p}$ is defined as (1.4).
\par  By uniqueness and decay estimates of the positive radial solutions (see \cite{M1993} and \cite{BL1983}), we state the following lemma.
\par \noindent\textbf{Lemma 2.5.} Let $\omega>0$, $d\geq 2$ and $1+\frac{4}{d} < p < \frac{d+2}{(d-2)^{+}}$. If (2.10) possesses a positive radial solution $Q_{\omega}(x)$, then $Q_{\omega}(x)$ has to be unique. Moreover $Q_{\omega}(x)$ has exponential decay property for some positive constant $C$:
\begin{align}
\rvert\nabla Q_{\omega}(x)\lvert+\rvert Q_{\omega}(x)\lvert\leq Ce^{-\frac{\sqrt{\omega}}{2}\rvert x\lvert}.
\end{align}


\section{Correlative Variational Framework}
\par  At first, we state a lemma on the Schwartz symmetric rearrangement (see \cite{LL2000}).
\par \noindent\textbf{Lemma 3.1.} For $\phi\in H^{1}$, let $\phi^{*}$ be the Schwartz symmetrization of $\phi$. Then for $d\geq 2$ and $1 < p < \frac{d+2}{(d-2)^{+}}$, one has that $\phi^{*}\geq 0$, a.e. in $\mathbb{R}^{d}$; $\phi^{*}(x)$ is radially symmetric; no-increasing and $\lim_{\rvert x\lvert\rightarrow\infty}\phi^{*}(x)\rightarrow 0$, $x\in \mathbb{R}^{d}$; in addition,
\begin{align}
\int\rvert\nabla \phi^{*}\lvert^{2}dx\leq\int\rvert\nabla\phi\lvert^{2}dx;\;\;\;\;\int\rvert \phi^{*}\lvert^{2}dx=\int\rvert \phi\lvert^{2}dx;
\end{align}
\begin{align}
\int\rvert \phi^{*}\lvert^{p+1}dx=\int\rvert\phi\lvert^{p+1}dx.
\end{align}

\par  Then we state  the radial compactness Lemma by Strauss \cite{S1977}.
\par \noindent\textbf{Lemma 3.2.} Put  $H_{r}^{1}(\mathbb{R}^{d})=\{\phi(x)\in H^{1},\;\;\phi(x)=\phi(\rvert x\lvert)\}$. Then for $d\geq 2$ and $1 < p < \frac{d+2}{(d-2)^{+}}$, the embedding $H_{r}^{1}(\mathbb{R}^{d})\hookrightarrow L^{p+1}(\mathbb{R}^{d})$ is compact.

\par  Now, we solve the following constrained variational problem.
\par \noindent\textbf{Theorem 3.3.} For $d\geq 2$, $1+\frac{4}{d} < p < \frac{d+2}{(d-2)^{+}}$ and $\int q^{2}dx<m<\infty$, where $q$ is the unique positive radial solution of (1.5), we set the constrained variational problem
\begin{align}
d_{m}:=inf_{\{u\in H^{1},\;\int \rvert u\lvert^{2}dx=m\}}E(u).
\end{align}
Then (3.3) possesses a positive minimizer $Q_{m} \in H^{1}$. Moreover there exists a unique $\omega_{m}>0$ such that $Q_{m}$ is the unique positive radial solution of (2.10) with $\omega=\omega_{m}$.
\\
\noindent\textbf{Proof.} Since
\begin{align}
\int q^{2}dx<\int \rvert u\lvert^{2}dx=m<\infty,
\end{align}
With $0<\varepsilon<\frac{2}{p+1}$, from (2.7) it follows that
\begin{align}
E(u) \geq -C(\varepsilon,\;d,\;p,\;\|u\|_{L^{2}}) > -\infty.
\end{align}
Then from (3.3)
\begin{align}
d_{m} > -\infty.
\end{align}
Taking $\mu =\frac{m}{\int q^{2}dx}$, from (3.4), $\mu > 1$. Let $u(x)=\mu q$. Then from (1.5) and (2.2)
\begin{align}
E(u)=\mu^{2} (1-\mu^{2})\int \rvert\nabla q\lvert^{2}dx+\frac{2}{p+1}\mu^{p+1}\int \rvert q\lvert^{p+1}dx.
\end{align}
Now, let $Q(x)=\lambda^{\frac{d}{2}} \mu q(\lambda x)$. Then there exists $0<\lambda<<1$ such that
$\int \rvert Q\lvert^{2}dx=m$ and $E(Q)<0$ since $\mu>1$. Therefore from (3.3), $d_{m} <0$. Combining with (3.6), we have that
\begin{align}
-\infty < d_{m} <0.
\end{align}
Then we can take a minimizing sequence $\{u_{n}\in H^{1}\}$ of the constrained variational problem (3.3) such that $\int \rvert u_{n}\lvert^{2}dx=m$ for $n\in \mathbb{N}$ and
\begin{align}
\lim_{n\rightarrow\infty}E(u_{n})=d_{m}.
\end{align}
Let $u^{*}_{n}$ be the Schwartz symmetrization of $u_{n}$ for any $n\in \mathbb{N}$. By (2.2), (3.3) and Lemma 3.1, it follows that $u^{*}_{n} \geq 0$, a.e. in $\mathbb{R}^{d}$, $\int \rvert u^{*}_{n}\lvert^{2}dx=m$ for $n\in \mathbb{N}$ and
\begin{align}
\liminf_{n\rightarrow\infty} E(u^{*}_{n})\leq d_{m}.
\end{align}
From (2.7), (3.4) and (3.6), it follows that a subsequence of $\{u^{*}_{n}\}$ still denoted by $\{u^{*}_{n}\}$ is a bounded sequence in $H^{1}$. Thus there exists a weak convergence subsequence of $\{u^{*}_{n}\}$ still denoted by $\{u^{*}_{n}\}$ such that as $n\rightarrow\infty$
\begin{align}
u^{*}_{n}\rightharpoonup v\;\;in\;\;H^{1},\;\;\;u^{*}_{n}(x)\rightarrow v(x)\;\;a.e.\;\; in\;\mathbb{R}^{d}.
\end{align}
From Lemma 3.2 it follows that
\begin{align}
u^{*}_{n}\rightarrow v\;\;in\;\;L^{2+\frac{4}{d}}(\mathbb{R}^{d}).
\end{align}
By weak lower semi-continuity, we obtain that
\begin{align}
\int\rvert v\lvert^{2}dx\leq m,\;\;\;\;E(v)\leq d_{m}.
\end{align}
(3.8) and (3.13) derive that $E(v) < 0$. It follows that $v\not\equiv 0$. Now let $v_{\lambda}=v(\lambda x)$
for $0<\lambda < \infty$, then
\begin{align}
\int\rvert v_{\lambda}\lvert^{2}dx=\lambda^{-d}\int\rvert v\lvert^{2}dx,
\end{align}
\begin{align}
E(v_{\lambda})=\lambda^{2-d}\int \rvert\nabla v\lvert^{2}dx-\lambda^{-d}\int(\frac{d}{d+2}\rvert v\lvert^{2+\frac{4}{d}}
-\frac{2}{p+1}\rvert v\lvert^{p+1})dx.
\end{align}
Since (3.13) and (3.14), there exists a $\lambda_{0}\in (0,\;1]$ such that
\begin{align}
\int\rvert v_{\lambda_{0}}\lvert^{2}dx=m.
\end{align}
From $E(v)<0$ and $\lambda_{0}\in (0,\;1]$, (3.15) derives that
\begin{align}
E(v_{\lambda_{0}}) \leq E(v) \leq d_{m}.
\end{align}
By (3.3) and (3.16), we have $E(v_{\lambda_{0}}) \geq d_{m}$. Then we deduce that $\lambda_{0} =1$ and
\begin{align}
\int\rvert v\lvert^{2}dx=m,\;\;E(v)=d_{m}.
\end{align}
Thus $v(x)$ is a minimizer of (3.3). Put $Q_{m}(x)=\rvert v(x)\lvert$. Then $Q_{m}(x) \geq 0$ is still a minimizer of (3.3). In terms of (3.3), there exists a unique Lagrange multiplier $\omega_{m}\in \mathbb{R}$ such that $Q_{m}$ has to satisfy the Euler-Lagrange equation
\begin{align}
\frac{d}{d\varepsilon}\lvert_{\varepsilon=0}[E(Q_{m}+\varepsilon \eta)+\omega_{m}\int\rvert Q_{m}+\varepsilon \eta\lvert^{2}dx-m\omega_{m}]=0,\;\;\;\;\eta\in C^{\infty}_{0}(\mathbb{R}^{d}).
\end{align}
It follows that $Q_{m}$ satisfies nonlinear elliptic equation (2.10) with $\omega=\omega_{m}$. Since $Q_{m}(x)=\rvert v(x)\lvert \geq 0$ a.e in $\mathbb{R}^{d}$, by the strong maximum principle we get that $Q_{m}(x)=\rvert v(x)\lvert > 0$ for $x\in \mathbb{R}^{d}$. Thus $Q_{m}(x)=\rvert v(x)\lvert$ is a positive minimizer of (3.3). By Theorem 2.3, the Lagrange multiplier $\omega_{m} > 0$. From Lemma 2.5, $Q_{m}$ is the unique positive radial solution of (2.10) with $\omega=\omega_{m}$.
\par  This completes the proof of Theorem 3.3.
\par \noindent\textbf{Theorem 3.4.} Let $q$ be the unique positive radial solution of (1.5) and $Q_{\omega}$ be the positive radial solutions of (2.10).
Define
\begin{align}
\mu_{q}=\{\omega\in \mathbb{R}\rvert\int Q^{2}_{\omega}dx > \int q^{2}dx\}.
\end{align}
Then $\mu_{q}$ is not empty. Moreover $0 \leq \omega_{q}=inf\mu_{q} < \omega_{p}$.
\par \noindent\textbf{Proof.} By Theorem 3.3, the Lagrange multiplier $\omega_{m} \in \mu_{q}$. It follows that $\mu_{q}$ is not empty. Then (1.9) and Lemma 2.4 deduce that $0 \leq \omega_{q}=inf\mu_{q} < \omega_{p}$.
\par  This proves Theorem 3.4.
\par \noindent\textbf{Theorem 3.5.} For $\int q^{2}dx < m <\infty$, (3.3) has a unique positive radial minimizer $Q_{m}\in H^{1}$. Moreover the set of all solutions of (3.3) is $S_{m}=\{e^{i\theta}Q_{m}(\cdot+y),\;\;\theta\in \mathbb{R},\;\;y\in \mathbb{R}^{d}\}$. In addition, for arbitrary $u\in S_{m}$, there exists a unique $\omega=\omega_{m}>0$ such that $\varphi(t,\;x)=e^{i\omega_{m} t}u(x)$ is a soliton of (1.1).
\par \noindent\textbf{Proof.} From Theorem 3.3, (3.3) has a positive radial minimizer $Q_{m} \in H^{1}$. Now suppose that $v\in H^{1}$ is an arbitrary solution of (3.3). Let $v=v^{1}+iv^{2}$, where $v^{1},\;\;v^{2}\in H^{1}$ are real-valued. Then $\widetilde{v}=\rvert v^{1}\lvert+i\rvert v^{2}\lvert$ is still a solution of (3.3). Thus there exists a unique $\omega_{m}>0$ such that $v$ and $\widetilde{v}$ satisfy (2.10). Then for $j=1,\;2$,
\begin{align}
\Delta v^{j}+\rvert v\lvert^{\frac{4}{d}}v^{j}-\rvert v\lvert^{p-1}v^{j}=\omega_{m} v^{j}\;\;in\;\;\mathbb{R}^{d},
\end{align}
\begin{align}
\Delta \rvert v^{j}\lvert+\rvert v\lvert^{\frac{4}{d}}\rvert v^{j}\lvert-\rvert v\lvert^{p-1}\rvert v^{j}\lvert=\omega_{m} \rvert v^{j}\lvert\;\;in\;\;\mathbb{R}^{d}.
\end{align}
This shows that $\omega_{m}$ is the first eigenvalue of the operator $\Delta+\rvert v\lvert^{\frac{4}{d}}-\rvert v\lvert^{p-1}$ acting over $H^{1}$ and thus $v^{1},\;v^{2},\;\rvert v^{1}\lvert$ and $\rvert v^{2}\lvert$ are all multiples of a positive normalized eigenfunction $v_{0}$ of $\Delta+\rvert v\lvert^{\frac{4}{d}}-\rvert v\lvert^{p-1}$, i.e.
\begin{align}
\Delta v_{0}+\rvert v\lvert^{\frac{4}{d}}v_{0}-\rvert v\lvert^{p-1}v_{0}=\omega_{m} v_{0} \;\;in\;\;\mathbb{R}^{d}
\end{align}
with
\begin{align}
v_{0}\in C^{2}(\mathbb{R}^{d})\cap H^{1},\;\;\;\;v_{0}>0\;\;in\;\;\mathbb{R}^{d}\;\;and\;\;\int\rvert v_{0}\lvert^{2}dx=m.
\end{align}
It is now obvious to deduce that: $v=e^{i\theta}v_{0}$ for some $\theta\in \mathbb{R}$ and that $v_{0}$ is still a solution of (3.3). In terms of Theorem 3.3, $v_{0}$ is the unique positive radial solution of (2.10) for the above $\omega=\omega_{m}$. It follows that $v_{0}=Q_{m}(\cdot+y)$ for some $y\in \mathbb{R}^{d}$. Thus $v=e^{i\theta}Q_{m}(\cdot+y)$ for some $\theta\in \mathbb{R}$. It is obvious that for any $\theta\in \mathbb{R}$ and $y\in \mathbb{R}^{d}$, $e^{i\theta}Q_{m}(\cdot+y)$ is also a solution of (3.3). Therefore
\begin{align}
S_{m}=\{e^{i\theta}Q_{m}(\cdot+y),\;\;\theta\in \mathbb{R},\;\;y\in \mathbb{R}^{d}\}
\end{align}
is the set of all solutions of (3.3). Moreover for arbitrary $u\in S_{m}$, there exists a unique $\omega=\omega_{m}>0$ such that $u$ is a solution of (2.10), which turns out that $\varphi(t,\;x)=e^{i\omega_{m} t}u(x)$ is a soliton solution of (1.1).
\par  This completes the proof of Theorem 3.5.
\par \noindent\textbf{Lemma 3.6.} For $d \geq 2$, $1+\frac{4}{d} < p <\frac{d+2}{(d-2)^{+}}$ and $u \in H^{1}\backslash \{0\}$, define the functional
\begin{align}
I(u)=2\int \rvert\nabla u\lvert^{2}dx-\frac{2d}{d+2}\int \rvert u\lvert^{2+\frac{4}{d}}dx+\frac{p-1}{p+1}d\int \rvert u\lvert^{p+1}dx.
\end{align}
For $\lambda > 0$, let $u_{\lambda}=\lambda^{\frac{d}{2}}u(\lambda x)$. Then to $\omega > 0$, we have that
\begin{align}
\frac{d}{d\lambda}[E(u_{\lambda})+\omega \int \rvert u_{\lambda}\lvert^{2}dx]=\frac{1}{\lambda}I(u_{\lambda}).
\end{align}
In addition $E(u_{\lambda})+\omega \int \rvert u_{\lambda}\lvert^{2}dx$ attains the minimum at $\lambda_{0}$ satisfying $I(u_{\lambda_{0}})=0$.
Moreover if $u$ is a solution of (2.10), one has that $I(u)=0$.
\par \noindent\textbf{Proof.}  By a direct calculation, it is shown that (3.27) is true. In addition
$E(u_{\lambda})+\omega \int \rvert u_{\lambda}\lvert^{2}dx$ attains the minimum at $\lambda_{0}$ satisfying $I(u_{\lambda_{0}})=0$.
Moreover if $u$ is a solution of (2.10), from (2.12), it follows that $I(u)=0$.
\par  This completes the proof of Lemma 3.6.
\par  Now we introduce and solve the following correlative constrained variational problem.
\par \noindent\textbf{Theorem 3.7.}  For $d \geq 2$, $1+\frac{4}{d} < p <\frac{d+2}{(d-2)^{+}}$, $\omega \in (\omega_{q}, \omega_{p})$ and $\int q^{2}dx < m < \infty$, we set the constrained variational problem
\begin{align}
d_{\omega}=inf_{\{u\in H^{1},\;\int q^{2}dx < \int \rvert u\lvert^{2}dx \leq m,\;I(u)=0\}}[E(u)+\omega\int\rvert u\lvert^{2}dx],
\end{align}
where $\omega_{p}$ is defined as (1.4) and $\omega_{q}$ is defined as (1.9).
Then (3.28) possesses a positive minimizer $Q_{\omega}\in H^{1}$. Moreover $Q_{\omega}$ is the unique positive radial solution of  (2.10).
\par \noindent\textbf{Proof.}  In the following we complete this proof by five steps.
\par  Step 1. $\{u\in H^{1},\;\int q^{2}dx < \int \rvert u\lvert^{2}dx \leq m,\;I(u)=0\}$ is not empty.
\par  Since $\int q^{2}dx < m < \infty$, from Theorem 3.3 it follows that there exists a positive minimizer $Q_{m}(x) \in H^{1}$ such that
$\int Q_{m}(x)^{2}dx = m$ and $Q_{m}(x)$ satisfies (2.10) with $\omega=\omega_{m}\in (\omega_{q},\;\omega_{p})$. By Lemma 3.6 it follows that
$I(Q_{m})=0$. Thus $Q_{m}\in\{u\in H^{1},\;\int q^{2}dx < \int \rvert u\lvert^{2}dx \leq m,\;I(u)=0\}$. Therefore the set $\{u\in H^{1},\;\int q^{2}dx < \int \rvert Q_{m}(x)\lvert^{2}dx \leq m,\;I(u)=0\}$ is not empty.
\par  Step 2. $d_{\omega} > -\infty$.
\par For $\int q^{2}dx < m < \infty$, take $u\in H^{1}$ satisfying $I(u)=0$ and $\int q^{2}dx < \int \rvert u\lvert^{2}dx \leq m$. For $\omega\in (\omega_{q},\;\omega_{p})$, from (2.7)
\begin{align}
E(u)+\omega \int \rvert u\lvert^{2}dx=&\int \rvert\nabla u\lvert^{2}dx-\frac{d}{d+2}\int\rvert u\lvert^{2+\frac{4}{d}}dx+\frac{2}{p+1}\rvert u\lvert^{p+1}dx+\omega\int \rvert u\lvert^{2}dx\notag\\
\geq &\int \rvert\nabla u\lvert^{2}dx-C(\varepsilon,d,p,\|u\|_{L^{2}})+(\frac{2}{p+1}-\varepsilon)\int\rvert u\lvert^{p+1}dx+\omega\int \rvert u\lvert^{2}dx.
\end{align}
Taking $0 < \varepsilon <\frac{2}{p+1}$, since $\int q^{2}dx < \int \rvert u\lvert^{2}dx \leq m$, by (3.29) we have that
\begin{align}
E(u)+\omega \int \rvert u\lvert^{2}dx \geq -C(\varepsilon,d,p,m)=constant>-\infty.
\end{align}
Therefore we deduce that $d_{\omega} > -\infty$.
\par  Step 3. Minimizing sequence is bounded in $H^{1}$.
\par  Let $\{u_{n}\}\subset H^{1}$ be a minimizing sequence of (3.28). Then for all $
n\in \mathbb{N}$, $I(u_{n})=0$, $\int q^{2}dx < \int \rvert u_{n}(x)\lvert^{2}dx \leq m$ and
\begin{align}
E(u_{n})+\omega\int\lvert u_{n}\rvert^{2}dx\rightarrow d_{\omega},  n\rightarrow\infty.
\end{align}
By (3.29) and (3.31), for $n$ large enough we have that
\begin{align}
\int\lvert \nabla u_{n}\rvert^{2}dx < d_{\omega}+1+C(\varepsilon,d,p,m).
\end{align}
Combining with $\int q^{2}dx < \int \rvert u_{n}(x)\lvert^{2}dx \leq m$, we deduce that $\{u_{n}\}$ is bounded in $H^{1}$.
\par  Step 4. Existence of minimizer.
\par  For a minimizing sequence $\{u_{n}\}$ of (3.28), let $u^{*}_{n}$ be the Schwartz symmetrization of $u_{n}$ for every $n\in \mathbb{N}$. By Lemma 3.1, $u^{*}_{n} \geq 0$, a.e. in $\mathbb{R}^{d}$ and $\{u^{*}_{n}\}$ is also bounded in $H^{1}$. In addition, for all $n\in \mathbb{N}$, $\int q^{2}dx < \int (u^{*}_{n})^{2}dx \leq m$,
\begin{align}
I(u^{*}_{n}) \leq 0,~~~E(u^{*}_{n})+\omega\int (u^{*}_{n})^{2}dx \leq E(u_{n})+\omega\int\lvert u_{n}\rvert^{2}dx.
\end{align}
Moreover there exists a weak convergence subsequence of $\{u^{*}_{n}\}$ stills denoted by $\{u^{*}_{n}\}$ such that
\begin{align}
u^{*}_{n}\rightharpoonup v~in~H^{1},~~u^{*}_{n}(x)\rightarrow v(x)~~a.e.~~in~~\mathbb{R}^{d}.
\end{align}
From Lemma 3.2, it follows that
\begin{align}
u^{*}_{n}\rightarrow~v~in~L^{2+\frac{4}{d}}(\mathbb{R}^{d}).
\end{align}
Applying weak lower semi-continuity, by (3.33), (3.34) and (3.35) we deduce that
\begin{align}
\int \rvert v\lvert^{2}dx \leq m,\;I(v)\leq 0,\;E(v)+\omega\int\lvert v\rvert^{2}dx \leq d_{\omega}.
\end{align}
\par  Now we claim that $v\not\equiv 0$.
\par  Otherwise, (3.34) derives that
\begin{align}
u^{*}_{n}\rightharpoonup~0~in~H^{1}, n\rightarrow \infty.
\end{align}
By (3.37), it deduces that
\begin{align}
\int (u^{*}_{n})^{2}dx \rightarrow 0, n\rightarrow \infty.
\end{align}
This is contradictory with $\int q^{2}dx < \int (u^{*}_{n})^{2}dx=\int \rvert u_{n}\lvert^{2}dx \leq m $ for any $n\in \mathbb{N}$.
\par  Now let $v_{\lambda}=\lambda^{\frac{d}{2}}v(\lambda x)$ for $\lambda > 0$. By $I(v) \leq 0$, there exists $\lambda_{0} > 0$ such that $I(v_{\lambda_{0}})=0$. From Lemma 3.6, we have that
\begin{align}
E(v_{\lambda_{0}})+\omega \int\rvert v_{\lambda_{0}}\lvert^{2}dx \leq E(v)+\omega \int\rvert v\lvert^{2}dx.
\end{align}
On the other hand, by $I(v_{\lambda_{0}})=0$, (3.26) and Lemma 2.2 derive that
\begin{align}
\int \rvert v_{\lambda_{0}}\lvert^{2}dx > \int q^{2}dx.
\end{align}
Combining with $\int \rvert v_{\lambda_{0}}\lvert^{2}dx=\int \rvert v\lvert^{2}dx$, it follows that
\begin{align}
\int  q^{2}dx < \int \rvert v_{\lambda_{0}}\lvert^{2}dx \leq m.
\end{align}
Therefore (3.28) derives that
\begin{align}
d_{\omega} \leq E(v_{\lambda_{0}})+\omega \int\rvert v_{\lambda_{0}}\lvert^{2}dx.
\end{align}
Thus (3.36), (3.39) and (3.42) yield that
\begin{align}
E(v)+\omega \int\rvert v\lvert^{2}dx=E(v_{\lambda_{0}})+\omega \int\rvert v_{\lambda_{0}}\lvert^{2}dx=d_{\omega}.
\end{align}
It follows that $\lambda_{0}=1$. Therefore we conclude that
\begin{align}
I(v)=0,\;\;E(v)+\omega \int\rvert v\lvert^{2}dx=d_{\omega}.
\end{align}
Combining with (3.41), we get that $v=v_{\lambda_{0}}$ with $\lambda_{0}=1$ satisfies $\int q^{2}dx < \int \rvert v\lvert^{2}dx \leq m$. Thus $v(x)$ is a minimizer of (3.28). Put
$Q_{\omega}(x)=\rvert v(x)\lvert$. Then it is clear that $Q_{\omega}(x)$ is still a minimizer of (3.28).
\par  Step 5. $Q_{\omega}(x)$ is the positive solution of (2.10).
\par  In terms of (3.28), there exists a unique $\Lambda \in \mathbb{R}$ such that $Q_{\omega}(x)$ satisfies the Euler-Lagrange equation for $\eta \in C^{\infty}_{0}(\mathbb{R}^{d})$
\begin{align}
\frac{d}{d\varepsilon}\lvert_{\epsilon =0}[E(Q_{\omega}+\varepsilon\eta)+\omega\int\rvert Q_{\omega}+\varepsilon\eta\lvert^{2}dx+\Lambda I(Q_{\omega}+\varepsilon\eta)]=0.
\end{align}
It follows that
\begin{align}
-\Delta Q_{\omega}-Q_{\omega}^{1+\frac{4}{d}}+Q_{\omega}^{p}+\omega Q_{\omega}+\Lambda[-2\Delta Q_{\omega}-2Q_{\omega}^{1+\frac{4}{d}}+\frac{d}{2}(p-1)Q_{\omega}^{p}]=0.
\end{align}
From (3.46), we have that
\begin{align}
\int [1+2\Lambda]\rvert\nabla Q_{\omega}\lvert^{2}+[-1-2\Lambda]Q_{\omega}^{2+\frac{4}{d}}+[1+\frac{d}{2}\Lambda (p-1)]Q_{\omega}^{p+1}+\omega Q_{\omega}^{2}dx = 0;
\end{align}
\begin{align}
\int\frac{d-2}{2}(1+2\Lambda)\rvert\nabla Q_{\omega}\lvert^{2}&-[1+2\Lambda]\frac{d^{2}}{d+4}Q_{\omega}^{2+\frac{4}{d}}\notag\\
&+[1+\frac{d}{2}\Lambda(p-1)]\frac{d}{p+1}Q_{\omega}^{p+1}+\frac{d}{2}\omega Q_{\omega}^{2}dx=0.
\end{align}
By $I(Q_{\omega}) = 0$, (3.47) and (3.48), we have that
\begin{align}
\Lambda\int-\frac{p-1}{p+1}[(p-1)d-4]\lvert Q_{\omega}\lvert^{p+1}dx=0.
\end{align}
Since $1+\frac{4}{d} < p < \frac{d+2}{(d-2)^{+}}$ and $Q_{\omega}(x)\not\equiv 0$, we have $\Lambda=0$.

 From (3.46), $Q_{\omega}(x) \geq 0$ is a solution of (2.10).
Using the strong maximum principle we get that $Q_{\omega}(x) > 0$ for $x\in \mathbb{R}^{d}$. In terms of Lemma 2.5, $Q_{\omega}(x)=\rvert v(x)\lvert$ is the unique positive radial solution of (2.10). It turns that $Q_{\omega}(x)$ is a positive minimizer of (3.28).
\par  This completes the proof of Theorem 3.7.

\par \noindent\textbf{Theorem 3.8.}  The variational problem (3.3) determines a one-to-one mapping between $m\in (\int q^{2}dx, \infty)$ and $\omega\in (\omega_{q}, \omega_{p})$. In detail, for $\omega\in (\omega_{q}, \omega_{p})$ and $m=\int Q^{2}_{\omega}dx$ with positive solution $Q_{\omega}(x)$ of (2.10) , one has that $\frac{dm}{d\omega}=\frac{d}{d\omega}\int Q^{2}_{\omega}dx \neq 0$.
\par \noindent\textbf{Proof.}  For arbitrary $m\in (\int q^{2}dx, \infty)$, in terms of Theorem 3.3, the variational problem (3.3) determines a positive $Q_{m}(x)\in H^{1}$ and a unique $\omega_{m}$ satisfying (2.10) with $\int Q^{2}_{m}dx=m$. By Lemma 2.4 and Theorem 3.4, $\omega_{m}\in (\omega_{q}, \omega_{p})$.
\par  Now suppose that there exists another $\omega'\in (\omega_{q}, \omega_{p})$ such that $\omega'\neq \omega_{m}$ and $\int Q^{2}_{\omega'}dx=m$ for the positive solution $Q_{\omega'}(x)$ of (2.10) with $\omega=\omega'$. By Lemma 2.4, $Q_{\omega'}(x) > 0$ is unique for (2.10) with $\omega=\omega'$. In addition, $\omega'\neq \omega_{m}$ leads that $Q_{\omega'}(x)\neq Q_{m}(x)$. From Theorem 3.5, $Q_{\omega'}(x)$ is not a minimizer of (3.3) since $Q_{m}$ is the unique positive minimizer of (3.3). On the other hand, according to Theorem 3.7, that $Q_{\omega'}(x)$ is the unique positive solution of (2.10) with $\omega=\omega'$ derives that $Q_{\omega'}(x)$ has to be the unique positive minimizer of the variational problem (3.28) with $\omega=\omega'$. We see that $Q_{m}$ satisfies (2.10) with $\omega=\omega_{m}$ and $Q_{\omega'}$ satisfies (2.10) with $\omega=\omega'$. By Lemma 3.6, it follows that $I(Q_{m})=0=I(Q_{\omega'})$.
\par  Summarizing the above facts, we get that
\begin{align}
\int Q^{2}_{\omega'}dx=\int Q^{2}_{m}dx=m;
\end{align}
\begin{align}
I(Q_{\omega'})=I(Q_{m})=0;
\end{align}
\begin{align}
Q_{m}\;is \;the\;minimizer\;of\;(3.3);
\end{align}
\begin{align}
Q_{\omega'}\;is \;the\;minimizer\;of\;(3.28)\;with\;\omega=\omega'.
\end{align}
Since $Q_{\omega'}$ is not a minimizer of (3.3), by (3.50) and (3.52), Theorem 3.3 derives
\begin{align}
E(Q_{m}) < E(Q_{\omega'}).
\end{align}
By (3.50), (3.51) and (3.53), Theorem 3.7 derives that
\begin{align}
E(Q_{\omega'})+\omega'\int Q^{2}_{\omega'}dx \leq E(Q_{m})+\omega'\int Q^{2}_{m}dx.
\end{align}
From (3.50), it is clear that (3.54) is contradictory with (3.55). Therefore it is necessary that $\omega_{m}=\omega'$.
\par  On the other hand, for arbitrary $\omega \in (\omega_{q}, \omega_{p})$, by Lemma 2.4 and Theorem 3.4, there exists a unique positive solution $Q_{\omega}(x)$ of (2.10) such that $\int Q^{2}_{\omega}dx=m\in (\int q^{2}dx, +\infty)$. Therefore the variational problem (3.3) determines a one-to-one mapping between $m \in (\int q^{2}dx, +\infty)$ and $\omega \in (\omega_{q}, \omega_{p})$. It turns that for $\omega \in (\omega_{q}, \omega_{p})$ and the positive solution $Q_{\omega}(x)$ of (2.10), we have that $\frac{dm}{d\omega}=\frac{d}{d\omega}\int Q^{2}_{\omega}dx\neq 0$.
\par  This completes the proof of Theorem 3.8.
\par \noindent\textbf{Theorem 3.9. }  Let $\mu$ be the set of all Lagrange multipliers corresponding to the all positive minimizers of (3.3).
Then $\mu = (\omega_{q},\omega_{p})$, where $\omega_{p}$ and $\omega_{q}$ are defined as (1.4) and (1.9) respectively. In addition $\mu_{q}=(\omega_{q},\omega_{p})$, where $\mu_{q}$ is defined as (3.20).
\par \noindent\textbf{Proof.}  Theorem 3.4 derives that $\mu\subset\mu_{q}\subset(\omega_{q},\omega_{p})$. Now suppose that $\omega \in (\omega_{q}, \omega_{p})$. By Theorem 3.8, this $\omega$ determines a unique $m\in(\int q^{2}dx,\infty)$, and this $m$ determines a unique $\omega_{m}\in \mu$. Then this $\omega_{m}$ can only be $\omega$, that is $\omega_{m}= \omega$. Thus $(\omega_{q},\omega_{p})\subset \mu$. Therefore $\mu =(\omega_{q},\omega_{p})$. It follows that $\mu =(\omega_{q},\omega_{p})=\mu_{q}$.
\par  This proves Theorem 3.9.
\par \noindent\textbf{Remark 3.10. }  Theorem 3.8 shows that for $\int q^{2}dx < m < \infty$, the normalized solution problem (2.10) with $\int \rvert u\lvert^{2}dx=m$ possesses a unique  $\omega_{m} \in (\omega_{q}, \omega_{p})$ such that $Q_{m}$ is the unique positive solution of (2.10) with $\omega=\omega_{m}$. It gives a positive answer that for (2.10), the mapping from the prescribed mass $m$ to the Lagrange multiplier, that is the soliton frequency $\omega$ is injective. Moreover the approach introduced here can be used to deal with more nonlinear Schr\"odinger equations.

\section{Orbital Stability of Big Solitons}
We prove a new compactness theorem as follows.
\par \noindent\textbf{Theorem 4.1. }  Let $\{u_{n}\}$ be a minimizing sequence of (3.3). Then there exists a subsequence of $\{u_{n}\}$, still denoted by
$\{u_{n}\}$ such that $u_{n}\rightarrow u$ in $H^{1}$.
\par \noindent\textbf{Proof.}  By Theorem 3.3, one has that
\begin{align}
\int \rvert u_{n}\lvert^{2} dx=m,\;\;\;n\in \mathbb{N};\;\;\;\lim_{n\rightarrow \infty} E(u_{n})=d_{m}.
\end{align}
From (2.7) and (4.1), it yields that $\{u_{n}\}$ is a bounded sequence in $H^{1}$. It follows that there exists a weak convergence subsequence of $\{u_{n}\}$ still denoted by $\{u_{n}\}$ such that as $n\rightarrow \infty$
\begin{align}
{u_{n}}\rightharpoonup u\;\;in\;\;H^{1},\;\;u_{n}\rightarrow u(x)\;\;a.e.\;\;in\;\;\mathbb{R}^{d}.
\end{align}
Let $u^{*}_{n}$ be the Schwartz symmetrization of $u_{n}$ for any $n\in \mathbb{N}$ and $u^{*}$ be the Schwartz symmetrization of $u$. Now we claim that
\begin{align}
u^{*}_{n}(x)\rightarrow u^{*}(x),\;a.e.\;\;in\;\;\mathbb{R}^{d}.
\end{align}
\par  In fact, we take $x\in \mathbb{R}^{d}$ such that $u_{n}(x)\rightarrow u(x)$ in $\mathbb{R}^{d}$. Then for arbitrary $\varepsilon >0$, there exists $N\in \mathbb{N}$ such that for $n > N$, $\rvert u_{n}(x)-u(x)\lvert < \varepsilon$. From \cite{LL2000},
\begin{align*}
\rvert u^{*}_{n}(x)-&u^{*}(x)\lvert=\rvert \int^{\infty}_{0}\chi_{\{\rvert u_{n}\lvert> \tau\}}^{*}(x)d\tau-\int^{\infty}_{0}\chi_{\{\rvert u\lvert> \tau\}}^{*}(x)d\tau\lvert\notag\\
=&\rvert\int^{\infty}_{0}[\chi_{\{\rvert u_{n}\lvert> \tau\}}^{*}(x)-\chi_{\{\rvert u\lvert> \tau\}}^{*}(x)]d\tau\lvert
\leq \int^{\infty}_{0}\chi_{\{\big\rvert \rvert u_{n}\lvert-\rvert u\lvert \big\lvert> \tau\}}^{*}(x)d\tau\notag\\
\leq& \int^{\infty}_{0}\chi_{\{\rvert u_{n}-u\lvert> \tau\}}^{*}(x)d\tau
\leq \int^{\varepsilon}_{0}1 d\tau=\varepsilon.
\end{align*}
It follows that  $u^{*}_{n}(x)\rightarrow u^{*}(x)$ in $\mathbb{R}^{d}$. Thus from (4.2), we get (4.3).
\par  By Lemma 3.1, $\{u^{*}_{n}\}$ is a bounded sequence in $H^{1}$. It follows that there exists a weak convergence subsequence of $\{u^{*}_{n}\}$ still denoted by $\{u^{*}_{n}\}$ such that
\begin{align}
u^{*}_{n}\rightharpoonup v\;\;in\;\;H^{1},\;\;u^{*}_{n}\rightarrow v(x),\;a.e.\;\;in\;\;\mathbb{R}^{d}.
\end{align}
By (4.3) and (4.4), we have that
\begin{align}
u^{*}(x)=v(x),\;\;a.e.\;in\;\;\mathbb{R}^{d}.
\end{align}
It is clear that $\{u^{*}_{n}\}$ is a minimizing sequence of (3.3), that is
\begin{align}
\int\rvert u_{n}^{*}\lvert^{2}dx=m.\;\;\;\;E(u_{n}^{*})\rightarrow d_{m}.
\end{align}
In terms of the proof of Theorem 3.3, by (4.4), (4.5), (4.6) and (3.18), we have that
\begin{align}
\int\rvert u^{*}\lvert^{2}dx=m.\;\;\;\;E(u^{*})=d_{m}.
\end{align}
From (4.4) and (4.5), Lemma 3.2 derives that
\begin{align}
\lim_{n\rightarrow \infty}\int\rvert u_{n}^{*}\lvert^{2+\frac{4}{d}}dx=\int\rvert u^{*}\lvert^{2+\frac{4}{d}}dx,
\end{align}
\begin{align}
\lim_{n\rightarrow \infty}\int\rvert u_{n}^{*}\lvert^{p+1}dx=\int\rvert u^{*}\lvert^{p+1}dx.
\end{align}
From (4.6), (4.7), (4.8) and (4.9), we deduce that
\begin{align}
\lim_{n\rightarrow \infty}\int\rvert\nabla u_{n}^{*}\lvert^{2}dx=\int\rvert\nabla u^{*}\lvert^{2}dx.
\end{align}
By Lemma 3.1, (4.1) and (4.7), we deduce that
\begin{align}
\lim_{n\rightarrow \infty}\int\rvert\nabla u_{n}\lvert^{2}dx=\int\rvert\nabla u^{*}\lvert^{2}dx.
\end{align}
From (4.2),
\begin{align}
\int\rvert\nabla u\lvert^{2}dx\leq\liminf_{n\rightarrow \infty}\int\rvert\nabla u_{n}\lvert^{2}dx.
\end{align}
Thus, (4.11), (4.12) and Lemma 3.1 conclude that
\begin{align}
\int\rvert\nabla u\lvert^{2}dx\leq\int\rvert\nabla u^{*}\lvert^{2}dx\leq \int \rvert\nabla u\lvert^{2}dx.
\end{align}
It follows that
\begin{align}
\int \rvert\nabla u\lvert^{2}dx=\int \rvert\nabla u^{*}\lvert^{2}dx.
\end{align}
(4.11) and (4.14) derive that
\begin{align}
\lim_{n\rightarrow \infty}\int\rvert\nabla u_{n}\lvert^{2}dx=\int\rvert\nabla u\lvert^{2}dx.
\end{align}
By (4.1), (4.6), Lemma 3.1 and (4.15), it follows that
\begin{align}
\|u_{n}\|_{H^{1}}\rightarrow \|u\|_{H^{1}},\;\;\;u\rightarrow \infty.
\end{align}
Then, by (4.2) and (4.16), we deduce that
\begin{align*}
u_{n}\rightarrow u\;\;in\;\;H^{1}.
\end{align*}
\par This completes the proof of Theorem 4.1.

\par \noindent\textbf{Theorem 4.2} The soliton $e^{i\omega t}u(x)$ in Theorem 3.5 holds the orbital stability, i.e.: for arbitrary $\varepsilon>0$, there exists $\delta>0$ such that for any $\varphi_{0}\in H^{1}$, if
\begin{align}
inf_{\{\theta\in \mathbb{R},\;\;y\in \mathbb{R}^{d}\}}\|\varphi_{0}(\cdot)-e^{i\theta}u(\cdot+y)\|_{H^{1}}<\delta,
\end{align}
then the solution $\varphi(t,\;x)$ of the Cauchy problem (1.1)-(2.1) satisfies
\begin{align}
inf_{\{\theta\in \mathbb{R},\;\;y\in \mathbb{R}^{d}\}}\|\varphi(t,\cdot)-e^{i\theta}u(\cdot+y)\|_{H^{1}}<\varepsilon,\;\;t\in \mathbb{R}.
\end{align}
\par \noindent\textbf{Proof.} By Theorem 3.5, it is clear that for arbitrary $u\in S_{m}$ one has that
\begin{align}
S_{m}=\{e^{i\theta}u(\cdot+y),\;\;\theta\in \mathbb{R},\;\;y\in \mathbb{R}^{d}\}.
\end{align}
In terms of Theorem 2.1, for $\varphi_{0}\in H^{1}$ satisfying $\int \rvert\varphi_{0}\lvert^{2}dx=m>\int q^{2}dx$, the Cauchy problem (1.1)-(2.1) possesses a unique global solution $\varphi(t,\;x)\in C(\mathbb{R,}\;H^{1})$ with the mass conservation
\begin{align*}
\int\rvert\varphi\lvert^{2}dx=\int\rvert\varphi_{0}\lvert^{2}dx=m>\int q^{2}dx.
\end{align*}
and energy conservation $E(\varphi)=E(\varphi_{0})$ for all $t\in \mathbb{R}$. Now arguing by contradiction.
\par  If the conclusion of Theorem 4.2 does not hold, then there exist $\varepsilon>0$, a sequence $(\varphi^{n}_{0})_{n\in \mathbb{N}^{+}}$ such that
\begin{align}
inf_{\{\theta\in \mathbb{R},\;\;y\in \mathbb{R}^{d}\}}\|\varphi^{n}_{0}-e^{i\theta}u(\cdot+y)\|_{H^{1}}<\frac{1}{n},
\end{align}
and a sequence $(t_{n})_{n\in \mathbb{N}^{+}}$ such that
\begin{align}
inf_{\{\theta\in \mathbb{R},\;\;y\in \mathbb{R}^{d}\}}\|\varphi_{n}(t_{n},\;\cdot)-e^{i\theta}u(\cdot+y)\|_{H^{1}}\geq\varepsilon,
\end{align}
where $\varphi_{n}$ denotes the solution of the Cauchy problem (1.1)-(2.1) with initial datum $\varphi^{n}_{0}$. From (4.20) we yield that
\begin{align}
\int\rvert\varphi^{n}_{0}\lvert^{2}dx\rightarrow\int\rvert u\lvert^{2}dx=m,
\end{align}
\begin{align}
E(\varphi^{n}_{0})\rightarrow E(u)=d_{m}.
\end{align}
Thus (4.22), (4.23) and the conservations of the mass and the energy derive that $(\varphi_{n}(t_{n},\;\cdot))_{n\in \mathbb{N}^{+}}$ is a minimizing sequence for the problem (3.3). Therefore Theorem 4.1, (4.22) and (4.23) derive that, there exists $\theta\in \mathbb{R}$ and $y\in \mathbb{R}^{d}$ such that
\begin{align}
\lim_{n\rightarrow \infty}\|\varphi_{n}(t_{n},\;\cdot)-e^{i\theta}u(\cdot+y)\|_{H^{1}}=0.
\end{align}
This is contradictory with (4.21). Thus Theorem 4.2 is proved.
\par \noindent\textbf{Theorem 4.3.} Let $\omega\in (0, \omega_{p})$ and $Q_{\omega}$ be a positive radial solution of (2.10). Then we have
\begin{align}
 \frac{d}{d\omega}E(Q_{\omega})=-\omega\frac{d}{d\omega}M(Q_{\omega})=-\omega\frac{d}{d\omega}\int Q^{2}_{\omega}dx.
\end{align}
\noindent\textbf{Proof.}  Since $Q_{\omega}$ is a positive solution of (2.10), it follows that
\begin{align}
\Delta Q_{\omega}+Q^{\frac{4}{d}+1}_{\omega}-Q^{p}_{\omega}-\omega Q_{\omega}=0,\;\;\;Q_{\omega}\in H^{1}.
\end{align}
From (2.2) and (4.26), we have
\begin{align}
\frac{d}{d\omega}E(Q_{\omega})
={} &\frac{d}{d\omega}(\int\rvert\nabla Q_{\omega}\lvert^{2}-\frac{1}{1+2/d}\rvert Q_{\omega}\lvert^{\frac{4}{d}+2}dx+\frac{2}{p+1}\rvert Q_{\omega}\lvert^{p+1}dx)\notag\\
={} &\int2\rvert\nabla Q_{\omega}\lvert\frac{d}{d\omega}\rvert\nabla Q_{\omega}\lvert-2Q^{\frac{4}{d}+1}_{\omega}\frac{d}{d\omega}Q_{\omega}+2Q^{p}_{\omega}\frac{d}{d\omega}Q_{\omega}dx\notag\\
={} &\int-2\Delta Q_{\omega}\frac{d}{d\omega}Q_{\omega}-2Q^{\frac{4}{d}+1}_{\omega}\frac{d}{d\omega}Q_{\omega}+2Q^{p}_{\omega}\frac{d}{d\omega}Q_{\omega}dx\notag\\
={} &\int-2\omega Q_{\omega}\frac{d}{d\omega}Q_{\omega}dx=-\omega\int\frac{d}{d\omega} Q^{2}_{\omega}dx=-\omega\frac{d}{d\omega}\int Q^{2}_{\omega}dx.
\end{align}
Noting (2.3), this proves (4.25).

Let $\omega\in (0, \omega_{p})$ and $Q_{\omega}(x)$ be the unique positive solution of (2.10). We set the scalar
\begin{align}
D(\omega)=E(Q_{\omega})+\omega M(Q_{\omega})
\end{align}
and the linearized operator of (4.26)
\begin{align}
H_{\omega}=-\Delta+\omega-(1+\frac{4}{d})Q^{\frac{4}{d}}_{\omega}+pQ^{p-1}_{\omega}.
\end{align}
It is clear that
\begin{align}
H_{\omega}=\frac{1}{2}E''(Q_{\omega})+\frac{1}{2}\omega M''(Q_{\omega}).
\end{align}
\par \noindent\textbf{Theorem 4.4.}  Let $\omega\in (0, \omega_{p})$ and $Q_{\omega}(x)$ be the unique positive solution of (2.10). Then the operator $H_{\omega}$ has one negative simple eigenvalue and has its kernel spanned by $iQ_{\omega}$. Moreover the positive spectrum of $H_{\omega}$ is bounded away from zero.
\par \noindent\textbf{Proof.} Since $\omega\in (0, \omega_{p})$, By Lemma 2.4, there exists a unique positive radial function $Q_{\omega}(x)$ satisfying (4.26). Now suppose that $\lambda\in \mathbb{R}$ satisfies $H_{\omega}Q_{\omega}=\lambda Q_{\omega}$, that is
\begin{align}
-\Delta Q_{\omega}+\omega Q_{\omega}-(1+\frac{4}{d})Q^{\frac{4}{d}+1}_{\omega}+pQ^{p}_{\omega}=\lambda Q_{\omega}.
\end{align}
From (4.26), it follows that
\begin{align}
-\frac{4}{d}Q^{\frac{4}{d}}_{\omega}+(p-1)Q^{p-1}_{\omega}=\lambda.
\end{align}
From Lemma 2.5, $lim_{\rvert x\lvert\rightarrow \infty} Q_{\omega}(x)=0$. Thus by $p > 1+\frac{4}{d}$, (4.32) yields that $\lambda < 0$. Moreover by (4.32), we can uniquely determine $\lambda$ as follows
\begin{align}
\lambda=\lambda_{-}=[\int-\frac{4}{d}Q^{\frac{4}{d}+2}_{\omega}+(p-1)Q^{p+1}_{\omega}dx]\big/\int Q^{2}_{\omega}dx.
\end{align}
Therefore we get that $H_{\omega}$ has one negative simple eigenvalue $\lambda_{-}$.  It follows that $H_{\omega}(iQ_{\omega})=\lambda_{-}(iQ_{\omega})$. From the uniqueness of $Q_{\omega}$ we get that the kernel is spanned by $iQ_{\omega}$.
\par  Now suppose that $\lambda > 0$ and $u\in H^{1}\backslash \{0\}$ satisfying $H_{\omega}u=\lambda u$, that is
\begin{align}
-\Delta u+\omega u-(1+\frac{4}{d})Q^{\frac{4}{d}}_{\omega}u+pQ^{p-1}_{\omega}u=\lambda u.
\end{align}
By Lemma 2.5,
\begin{align}
-(1+\frac{4}{d})Q^{\frac{4}{d}}_{\omega}+pQ^{p-1}_{\omega} : = g(x)=o(\rvert x\lvert^{-1}).
\end{align}
From Kato \cite{K1959}, $-\Delta+g(x)$ has no positive eigenvalues. Thus (4.34) derives that $\lambda \leq \omega$. By Weyl's theorem on the essential spectrum, the rest of the spectrum of $H_{\omega}$ is bounded away from zero (see \cite{RS1978}).
\par  This proves Theorem 4.4.
\par  By Theorem 4.4, $H_{\omega}$ with $T'(0)=i$ satisfies Assumption 3 in \cite{GSS1987} for $\omega\in (0,\omega_{p})$. With $J=-i$, $X=H^{1}$ and $E$ defined as (2.2), by Theorem 2.1 and Lemma 2.4, (1.1) satisfies Assumption 1 and 2 in \cite{GSS1987} for $\omega\in (0,\omega_{p})$. Thus we can use Theorem 4.7 in \cite{GSS1987} and get the following lemma.
\par \noindent\textbf{Lemma 4.5.} Let $\omega\in (0, \omega_{p})$ and $Q_{\omega}(x)$ be the unique positive solution of (2.10). If $D''(\omega)=\frac{d^{2}}{d\omega^{2}}D(\omega) < 0$, then the soliton $e^{i\omega t}Q_{\omega}(x)$ of (1.1) is unstable.
\par  Therefore we get the following theorem.
\par \noindent\textbf{Theorem 4.6.} Let $\omega\in (\omega_{q}, \omega_{p})$ and $Q_{\omega}(x)$ be the unique positive solution of (2.10). Then we have that
\begin{align}
\frac{d}{d\omega}\int Q_{\omega}^{2}dx>0.
\end{align}
\par \noindent\textbf{Proof.} From Theorem 4.3 and (4.28), we have that
\begin{align}
D''(\omega)=\frac{d}{d\omega}\int Q_{\omega}^{2}dx.
\end{align}
Since $\omega\in (\omega_{q}, \omega_{p})$, from Theorem 3.9 it follows that $\omega\in \mu$. In terms of Theorem 4.2, the soliton $e^{i\omega t}Q_{\omega}(x)$ holds the orbital stability. By Lemma 4.5, we deduce that $D''(\omega) \geq 0$. From (4.37) it follows that $\frac{d}{d\omega}\int Q_{\omega}^{2}dx \geq 0$. Set $m(\omega)=\int Q^{2}_{\omega}dx$. From $\omega\in (\omega_{q}, \omega_{p})$, Theorem 3.8 deduces that $\frac{d}{d\omega}\int Q_{\omega}^{2}dx = \frac{dm}{d\omega} \neq 0$. Therefore we get that $\frac{d}{d\omega}\int Q_{\omega}^{2}dx > 0$.
\par  This proves Theorem 4.6.
\par \noindent\textbf{Proof of Theorem A.} Theorem 3.8 and Theorem 4.6 deduce that Theorem A is true.
\par \noindent\textbf{Proof of Theorem B.} By Theorem 3.3 and Theorem 4.2, it follows that Theorem B is true. On the other hand, Theorem 4.4 and Theorem 4.6 also deduce that Theorem B is true.
\section{Construction of Multi-Solitons }
It is clear that (1.1) admits the following symmetries.
\par  Time-space translation invariance: if $\varphi(t,x)$ satisfies (1.1), then for any $(t_{0},\;x_{0})\in \mathbb{R}\times \mathbb{R}^{d}$
\begin{align}
\psi(t,x)=\varphi(t-t_{0},x-x_{0})
\end{align}
also satisfies (1.1).
\par  Phase invariance: if $\varphi(t,x)$ satisfies (1.1), then for any $\gamma_{0}\in \mathbb{R}$,
\begin{align}
\psi(t,x)=\varphi(t,x)e^{i\gamma_{0}}
\end{align}
also satisfies (1.1).
\par  Galilean invariance: if $\varphi(t,x)$ satisfies (1.1), then for any $v_{0}\in \mathbb{R}^{d}$,
\begin{align}
\psi(t,x)=\varphi(t,x-v_{0}t)e^{i(\frac{1}{2}v_{0}x-\frac{1}{4}\rvert v_{0}\lvert^{2}t)}
\end{align}
also satisfies (1.1).
\par Let $d\geq 2$, $1+\frac{4}{d} < p < \frac{d+2}{(d-2)^{+}}$, $\omega_{p}$ defined by (1.4) and $\omega_{q}$ defined by (1.9). For $K\geq 2$ and $k=1,\;2,\;\cdot\cdot\cdot,\;K$, we take $\omega^{0}_{k}\in (\omega_{q},\;\omega_{p}),\;\;\gamma^{0}_{k}\in \mathbb{R},\;\;x^{0}_{k}\in \mathbb{R}^{d}$ and $v_{k}\in \mathbb{R}^{d}$ with $v_{k}\neq v_{k'}$ to $k\neq k'$. By Theorem B,
\begin{align}
e^{i\omega^{0}_{k}t}Q_{\omega^{0}_{k}}(x),\;\;k=1,\;2,\;\cdot\cdot\cdot,\;K
\end{align}
are the stable solitons of (1.1). Then in terms of the above symmetries, for $k=1,\;2,\;\cdot\cdot\cdot,\;K$,
\begin{eqnarray}
R_{k}(t,x)=Q_{\omega_{k}^{0}}(x-x^{0}_{k}-v_{k}t)e^{i(\frac{1}{2}v_{k}x-\frac{1}{4}\rvert v_{k}\lvert^{2}t+\omega_{k}^{0}t+\gamma_{k}^{0})},\;(t,x)\in\ \mathbb{R}\times \mathbb{R}^{d}
\end{eqnarray}
are also the solitons of (1.1).

Now we set
\begin{align}
R(t)=\sum_{k=1}^{K}R_{k}(t,\;\cdot),\;\;t\in \mathbb{R}.
\end{align}
\par  At first we state the main theorem in this section, that is the existence of multi-solitons of (1.1).
\par \noindent\textbf{Theorem 5.1.} Let $d\geq 2$, $1+\frac{4}{d} < p < \frac{d+2}{(d-2)^{+}}$.
For $K\geq 2$, $k= {1,\cdot\cdot\cdot,K}$, taking $\omega^{0}_{k}\in (\omega_{q},\;\omega_{p})$, $\gamma_{k}^{0}\in \mathbb{R},\;x_{k}^{0}\in \mathbb{R}^{d}, \;v_{k}\in \mathbb{R}^{d}$ with $v_{k}\neq v_{k^{'}} \;\;to \;\;k\neq k'$ and
\begin{align}
R_{k}(t,x)=Q_{\omega_{k}^{0}}(x-x^{0}_{k}-v_{k}t)e^{i(\frac{1}{2}v_{k}x-\frac{1}{4}\rvert v_{k}\lvert^{2}t+\omega_{k}^{0}t+\gamma_{k}^{0})}
\end{align}
with $(t,x)\in\ \mathbb{R}\times \mathbb{R}^{d}$, there exists a solution $\varphi(t,\;x)$ of (1.1) such that
\begin{align}
\forall t\geq 0,\;\|\varphi(t)-\sum_{k=1}^{K}R_{k}(t)\|_{H^{1}}\leq Ce^{-\theta_{0}t}
\end{align}
for some $\theta_{0}>0$ and $C>0$.
\par \noindent\textbf{Proof of Theorem C.}  Theorem 5.1 directly deduces that Theorem C is true.
\par  In the following, according to  Martel, Merle and Tsai's way (see \cite{MM2006} and \cite{MMT2006}), we prove Theorem 5.1.
\par  Let $T_{n}>0$, $n=1,2,\cdot\cdot\cdot$ and $\lim_{n\rightarrow\infty}T_{n}=+\infty$. For $n=1,2,\cdot\cdot\cdot$, by Theorem 2.1, we can let $\varphi_{n}$ be the unique global solution in $H^{1}$ for the Cauchy problem
\begin{equation}\label{AH1}
\left\{
\begin{split}
&i\partial_{t}\varphi_{n}+\Delta \varphi_{n}+\rvert\varphi_{n}\lvert^{\frac{4}{d}}\varphi_{n}-\rvert\varphi_{n}\lvert^{p-1}\varphi_{n}=0,&\qquad  (t,x)\in\ \mathbb{R} \times \mathbb{R}^{d},\\
&\varphi_{n}(T_{n},x)=R(T_{n}).
\end{split}
\right.
\end{equation}
\par  We first state the following claim.
\par \noindent\textbf{Claim 5.2.} (Claim 1 in \cite{MM2006}) Let $\{v_{k}\}$, $k=1,\cdot\cdot\cdot,K$ be $K$ vectors of $\mathbb{R}^{d}$ such that for any $k\neq k',v_{k}\neq v_{k'}$.
Then, there exists an orthonormal basis $(e_{1},...,e_{d})$ of $\mathbb{R}^{d}$ such that for any  $k\neq k',(v_{k},e_{1})\neq (v_{k'},e_{1})$.

\par Without any restriction, we can assume that the direction $e_{1}$ given by Claim 5.2 is $x_{1}$, since (1.1) is invariant by rotation. Therefore, we may assume that for any $k\neq k',v_{k,1}\neq v_{k',1}$. We suppose in fact that
\begin{eqnarray}
v_{1,1}<v_{2,1}<\cdot\cdot\cdot<v_{K,1}.
\end{eqnarray}
Since $\omega^{0}_{k}\in (\omega_{q},\;\omega_{p})$ with $k=1,\cdot\cdot\cdot,K$, we can set $\theta_{0}>0$ such that
\begin{equation}\label{AH4}
\sqrt{\theta_{0}}=\frac{1}{16}min\{v_{2,1}-v_{1,1},\cdot\cdot\cdot,v_{K,1}-v_{K-1,1},\;\;\sqrt{\omega_{1}^{0}},\cdot\cdot\cdot,\sqrt{\omega_{K}^{0}}\}.
\end{equation}
\par Now we state the following uniform estimates about the sequence $\{\varphi_{n}\}$ in (5.9), which is the key point of the proof of Theorem 5.1.
\par \noindent\textbf{Proposition 5.3.} There exist $T_{0}>0,C_{0}>0$ and $\theta_{0}>0$ such that for all $n \geq 1$,
\begin{equation}\label{AH4}
\forall t\in [T_{0},T_{n}],\;\;\;\|\varphi_{n}(t)-R(t)\|_{H^{1}}\leq C_{0}e^{-\theta_{0}t}.
\end{equation}
\par  By Proposition 5.3, the sequence $\{\varphi_{n}\}$ has the following global bounded property.
\par \noindent\textbf{Lemma 5.4.} There exists a constant number $C > 0$ such that for any $t\in[T_{0},T_{n}]$ and all $n \geq 1$,
$\|\varphi_{n}(t)\|_{H^{1}}\leq C$.
\par \noindent\textbf{Proof.} It is clear that $\|R(t)\|_{H^{1}}$ is uniformly bounded by (5.6). From Proposition 5.3, we obtain the consequence directly.
\par  Now we claim a strong compactness result in $L^{2}$.
\par \noindent\textbf{Claim 5.5.} (Claim 9 in \cite{CC2011}) Take $\epsilon_{0}>0$. There exists $K_{0}=K_{0}(\epsilon_{0})>0$ such that for all $n$ large enough, we have
\begin{align}
\int_{\rvert x\lvert>K_{0}}\rvert\varphi_{n}(T_{0},x)\lvert^{2}dx \leq \epsilon_{0}.
\end{align}

\par \noindent\textbf{Lemma 5.6.} There exists $\psi_{0}\in H^{1}$ such that up to a subsequence
\begin{align}
\varphi_{n}(T_{0})\rightarrow \psi_{0}\;\;in\;H^{s}(\mathbb{R}^{d})\;as\;n\rightarrow +\infty
\end{align}
for any $0\leq s<1$.
\par \noindent\textbf{Proof.} By Lemma 5.4, there exists $\psi_{0}\in H^{1}$ such that up to a subsequence,
\begin{align}
\varphi_{n}(T_{0})\rightharpoonup \psi_{0}\;\;\;in\;H^{1}\;as\;n\rightarrow +\infty.
\end{align}
From Claim 5.5, it follows that
\begin{align}
\varphi_{n}(T_{0})\rightarrow \psi_{0}\;\;\;in\;L^{2}_{loc}(\mathbb{R}^{d})\;as\;n\rightarrow +\infty.
\end{align}
We conclude that
\begin{align}
\varphi_{n}(T_{0})\rightarrow \psi_{0}\;\;\;in\;L^{2}\;as\;n\rightarrow +\infty.
\end{align}
By interpolation we get (5.14).
\par  This proves Lemma 5.6.
\par  In terms of the uniform estimates Proposition 5.3 and the compactness result Lemma 5.6, we prove Theorem 5.1.
\par \noindent\textbf{Proof of Theorem 5.1.}  Let $\psi_{0}$ be given by Lemma 5.6.
From $1+\frac{4}{d} < p < \frac{d+2}{(d-2)^{+}}$, there exists $0<\sigma<1$ such that $1+\frac{4}{d} < p <1+\frac{4}{d-2\sigma}$ and
\begin{align}
\rvert(\rvert z_{1}\lvert^{\frac{4}{d}}z_{1}-\rvert z_{1}\lvert^{p-1}z_{1})&-(\rvert z_{2}\lvert^{\frac{4}{d}}z_{2}-\rvert z_{2}\lvert^{p-1}z_{2})\lvert\notag\\
&\leq C(1+\rvert z_{1}\lvert^{p-1}+\rvert z_{2}\lvert^{p-1})\rvert z_{1}-z_{2}\lvert
\end{align}
for all $z_{1},\;z_{2}\in \mathbb{C}$.
This implies that the Cauchy problem of (1.1) with $\varphi(T_{0},\;x)=\psi_{0}$ is well-posed in $H^{\sigma}(\mathbb{R}^{d})$ (see Theorem 5.1.1 in \cite{ C2003}, also refer to \cite{CW1990} and \cite{T1987}). Then we let $\varphi(t,\;x)\in C([T_{0},\;T],\;H^{\sigma}(\mathbb{R}^{d}))$ be the corresponding maximal solution of (1.1) with $\varphi(T_{0},\;x)=\psi_{0}$.  Combining with Lemma 5.6, we can obtain
\begin{align*}
\varphi_{n}(t)\rightarrow \varphi(t)\;\;in\;\;H^{\sigma}(\mathbb{R}^{d})\;\;as\;\;n\rightarrow +\infty
\end{align*}
for any $t\in [T_{0},\;T)$. By boundedness of $\varphi_{n}(t)\;\;in\;\;H^{1}$, we also have
\begin{align*}
\varphi_{n}(t)\rightharpoonup \varphi(t)\;\;in\;\;H^{1}\;\;as\;\;n\rightarrow +\infty
\end{align*}
for any $t\in [T_{0},\;T)$. By Proposition 5.3, for any $t\in [T_{0},\;T)$, we have
\begin{align}
\|\varphi(t)-R(t)\|_{H^{1}}\leq \liminf_{n\rightarrow \infty}\|\varphi_{n}(t)-R(t)\|_{H^{1}}\leq C_{0}e^{-\theta_{0}t}.
\end{align}
In particular, since $R(t)$ is bounded in $H^{1}$, there exists $C>0$ such that for any $t\in [T_{0},\;T)$ we have
\begin{align}
\|\varphi(t)\|_{H^{1}}\leq C_{0}e^{-\theta_{0}t}+\|-R(t)\|_{H^{1}}\leq C.
\end{align}
Recall that, by the blow up alternative (see \cite{C2003}), either $T=+\infty$ or $T<+\infty$ and
$lim_{t\rightarrow T}\|\varphi(t)\|_{H^{1}}=+\infty$. Therefore (5.20) implies that $T=+\infty$. From (5.19) we infer that for all $t\in [T_{0},\;+\infty)$ we have
\begin{align*}
\|\varphi(t)-R(t)\|_{H^{1}}\leq C_{0}e^{-\theta_{0}t}.
\end{align*}
\par  This completes the proof of Theorem 5.1.
\par The proof of the uniform estimates Proposition 5.3 relies on a bootstrap argument. We first state the following bootstrap result.
\par \noindent\textbf{Proposition 5.7.} There exist $A_{0}>0,\theta_{0}>0,T_{0}>0$ and $N_{0}>0$ such that for all $n\geq N_{0}$ and $t^{*}\in [T_{0},T_{n}]$, if
\begin{align}
\forall t\in [t^{*},T_{n}],\;\;\;\|\varphi_{n}(t)-R(t)\|_{H^{1}}\leq A_{0}e^{-\theta_{0}t},
\end{align}
then
\begin{equation}\label{AH4}
\forall t \in [t^{*},T_{n}],\;\;\;\|\varphi_{n}(t)-R(t)\|_{H^{1}}\leq \frac{A_{0}}{2}e^{-\theta_{0}t}.
\end{equation}
By  Proposition 5.7, we deduce the uniform estimates Proposition 5.3.
\par \noindent\textbf{Proof of Proposition 5.3.}(Proposition 1 in \cite{MM2006})
Let $t^{*}$ be the minimal time such that (5.21) holds:
\begin{align*}
t^{*} = min\{\tau\in[T_{0},T_{n}];\;(5.21)\;holds\;for\;all\;t\in[\tau,T_{n}]\}.
\end{align*}
We prove by contradiction that $t^{*}=T_{0}$. Indeed, assume that $t^{*}>T_{0}$. Then
\begin{align*}
\|\varphi_{n}(t^{*})-R(t^{*})\|_{H^{1}}\leq A_{0}e^{-\theta_{0}t},
\end{align*}
and by Proposition 5.7 we can improve this estimate in
\begin{align*}
\|\varphi_{n}(t^{*})-R(t^{*})\|_{H^{1}}\leq \frac{A_{0}}{2}e^{-\theta_{0}t}.
\end{align*}
Hence, by continuity of $\varphi_{n}(t)$ in $H^{1}$, there exists $T_{0}\leq t^{**}<t^{*}$ such that (5.21) holds for all $t\in [t^{**},\;t^{*}]$. This contradicts the minimality of $t^{*}$.
\par  This completes the proof of Proposition 5.3.
\par Now we begin to prove Proposition 5.7 by a series of lemmas.
\par  For $k=1,\;\cdot\cdot\cdot,\;K$, let $\omega_{k}(t)\in (\omega_{q},\;\omega_{p})$ and $Q_{\omega_{k}(t)}(x)$ be the unique positive radial solutions of (2.10). To $ x^{0}_{k},\;x_{k}(t),\;v_{k}\in \mathbb{R}^{d}$ and $\gamma_{k}(t)\in \mathbb{R},\;k=1,\cdot\cdot\cdot,K$, we set $\widetilde{x}_{k}(t)=x^{0}_{k}+v_{k}t+x_{k}(t),\;\delta_{k}(t)=-\frac{1}{4}\rvert v_{k}\lvert^{2}t+\omega^{0}_{k}t+\gamma_{k}(t)$,
\begin{align}
\widetilde{R}_{k}(t)=Q_{\omega_{k}(t)}(\cdot-\widetilde{x}_{k}(t))e^{i(\frac{1}{2}v_{k}x+\delta_{k}(t))},
\end{align}
\begin{eqnarray}
\widetilde{R}(t)=\sum^{K}_{k=1}\widetilde{R}_{k}(t)\;\;\;\;and\;\;\;\;\;\;\varepsilon(t,\cdot)=\varphi_{n}(t,\cdot)-\widetilde{R}(t).
\end{eqnarray}
\par \noindent\textbf{Lemma 5.8.}(Lemma 3 in \cite{MM2006})  There exists $C_{1}>0$ such that if $T_{0}$ is large enough, then there exists a unique $C^{1}$ function $(\omega_{k},x_{k},\gamma_{k}):[t^{*},T_{n}]\rightarrow(\omega_{q},\omega_{p})\times \mathbb{R}^{d}\times \mathbb{R}$, for any $k=1,2,\cdot\cdot\cdot,K$ such that
\begin{equation}\label{AH4}
Re\int \widetilde{R}_{k}(t)\overline{\varepsilon}(t)dx=Im\int \widetilde{R}_{k}(t)\overline{\varepsilon}(t)dx=0,\;\;Re\int \nabla \widetilde{R}_{k}(t)\overline{\varepsilon}(t)dx=0,
\end{equation}
\begin{equation}\label{AH4}
\|\varepsilon(t)\|_{H^{1}}+\sum_{k=1}^{K}\rvert \omega_{k}(t)-\omega_{k}^{0}\lvert\leq C_{1}A_{0}e^{-\theta_{0}t},
\end{equation}
\begin{equation}\label{AH4}
\rvert\dot{\omega}_{k}(t)\lvert^{2}+\rvert\dot{x}_{k}(t)\lvert^{2}+\rvert\dot{\gamma}_{k}(t)-(\omega_{k}(t)-\omega_{k}^{0})\lvert^{2}\leq C_{1}\|\varepsilon(t)\|^{2}_{H^{1}}+C_{1}e^{-2\theta_{0}t}.
\end{equation}
\noindent\textbf{Claim 5.9.} (Claim 2 in \cite{MM2006}) Let $z(t)\in H^{1}$ be a solution of (1.1). Let $h:x_{1}\in \mathbb{R}\mapsto h(x_{1})$ be a $C^{3}$ real-valued function of one variable such that $h$, $h'$ and $h'''$ are bounded. Then, for all $t\in \mathbb{R}$
\begin{eqnarray}
\frac{1}{2}\frac{d}{dt}\int\rvert z\lvert^{2}h(x_{1})dx=Im\int\partial_{x_{1}}z\bar zh'(x_{1})dx,
\end{eqnarray}
\begin{align}
\frac{1}{2}\frac{d}{dt}Im\int\partial_{x_{1}}z\bar zh(x_{1})dx=&\int\rvert\partial_{x_{1}}z\lvert^{2}h'(x_{1})dx-\frac{1}{4}\int\rvert z\lvert^{2}h'''(x_{1})dx\notag\\
-&\frac{1}{d+2}\int\rvert z\lvert^{2+\frac{4}{d}}h'(x_{1})dx+\frac{p-1}{2(p+1)}\int\rvert z\lvert^{p+1}h'(x_{1})dx,
\end{align}
and for $j=2,\;\cdot\cdot\cdot,\;d$,
\begin{eqnarray}
\frac{1}{2}\frac{d}{dt}Im\int\partial_{x_{j}}z\bar zh(x_{1})dx=Re\int\partial_{x_{j}}z\partial_{x_{1}}\bar zh'(x_{1})dx.
\end{eqnarray}
\par    Since $\varphi_{n}(T_{n})=R(T_{n})$ and at time $t=T_{n}$  the decomposition in (5.24) is unique, it follows that
\begin{equation}\label{AH4}
\varepsilon(T_{n})\equiv 0,\;\;\widetilde{R}(T_{n})\equiv R(T_{n}),\;\;\omega_{k}(T_{n})=\omega_{k}^{0},\;\;x_{k}(T_{n})=0,\;\;\gamma_{k}(T_{n})=\gamma_{k}^{0}.
\end{equation}
\par  Let $Y(s)$ be a $C^{3}$ function such that
\begin{equation}\label{AH4}
0\leq Y\leq 1,\;Y'\geq 0\;on\;\mathbb{R};\;\;\;Y(s)=0\;\;for\;s\leq -1;\;\;\;Y(s)=1\;\;for\;s>1;\;\;\;
\end{equation}
and satisfying for some constant $C>0$,
\begin{eqnarray*}
(Y'(s))^{2}\leq C Y(s),\;\;\;(Y''(s))^{2}\leq C Y'(s)\;\;\;for\;all\;s\in \mathbb{R}.
\end{eqnarray*}
Then consider
\begin{eqnarray*}
Y(s)=\frac{1}{16}(1+s)^{4}\;\;for\;\;s \geq -1\;\;close\;\;to\;\;-1,
\end{eqnarray*}
and
\begin{eqnarray*}
Y(s)=1-\frac{1}{16}(1-s)^{4}\;\;for\;\;s \leq 1\;\;close\;\;to\;\;1.
\end{eqnarray*}
For all $k=2,\cdot\cdot\cdot,K$, let
\begin{eqnarray*}
\sigma_{k}=\frac{1}{2}(v_{k-1,1}+v_{k,1}).
\end{eqnarray*}
For $L>0$ large enough to be fixed later, for any $k=2,\cdot\cdot\cdot,K-1$, let
\begin{eqnarray}
y_{k}(t,x)=Y(\frac{x_{1}-\sigma_{k}{t}}{L})-Y(\frac{x_{1}-\sigma_{k+1}{t}}{L}),
\end{eqnarray}
\begin{equation}
y_{1}(t,x)=1-Y(\frac{x_{1}-\sigma_{2}{t}}{L}),\;\;y_{K}(t,x)=Y(\frac{x_{1}-\sigma_{K}{t}}{L}).
\end{equation}
Finally, set for all $k=1,\;\cdot\cdot\cdot,\;K$:
\begin{equation}
I_{k}(t)=\int\rvert\varphi_{n}(t,x)\lvert^{2}y_{k}(t,x)dx,M_{k}(t)=Im \int \nabla \varphi_{n}(t,x)\bar \varphi_{n}(t,x)y_{k}(t,x)dx.
\end{equation}
The quantities $I_{k}(t)$ and $M_{k}(t)$ are local versions of the $L^{2}$ norm and momentum. Ordering the $v_{k,\;1}$ as in (5.10) was useful to split the various solitons using only the coordinate $x_{1}$.
\par  From (5.31) to (5.35), (5.21) and Claim 5.9 deduce the following lemma.
\par \noindent\textbf{Lemma 5.10.} (Lemma 4 in \cite{MM2006}) Let $L>0$. There exists $C>0$ such that if $L$ and $T_{0}$ are large enough, then for all $k=2,\;\cdot\cdot\cdot,\;K$, $t\in[t^{*},T_{n}]$, we have
\begin{eqnarray}
\rvert I_{k}(T_{n})-I_{k}(t)\lvert+\rvert M_{k}(T_{n})-M_{k}(t)\lvert\leq \frac{CA_{0}^{2}}{L}e^{-2\theta_{0}t},\;\;\;k=2,\cdot\cdot\cdot,K.
\end{eqnarray}

\par \noindent\textbf{Claim 5.11.} (Claim 3 in \cite{MM2006}) There exists $C>0$ such that for any $t\in[t^{*},\;T_{n}]$,
\begin{equation}
\rvert\omega_{k}(t)-\omega^{0}_{k}\lvert\leq C\|\varepsilon(t)\|^{2}_{L^{2}}+C(\frac{A_{0}^{2}}{L}+1)e^{-2\theta_{0}t}.
\end{equation}

\noindent\textbf{Proof.} From (5.24) and (5.35), we have
\begin{eqnarray*}
I_{k}(t)=\int\rvert\tilde{R}(t)\lvert^{2}y_{k}(t)dx+2Re \int\tilde{R}(t)\overline{\varepsilon}(t)y_{k}(t)dx+\int\rvert\varepsilon(t)\lvert^{2}y_{k}(t)dx.
\end{eqnarray*}
By the exponential decay of each $Q_{\omega_{k}(t)}$, the orthogonality $\int\tilde{R}_{k}(t)\bar{\varepsilon}(t)dx=0$ and the property of support of $y_{k}$, we have
\begin{eqnarray*}
I_{k}(t)=\int\rvert\varphi_{n}(t)\lvert^{2}y_{k}(t)dx=\int Q^{2}_{\omega_{k}(t)}dx+\int \rvert\varepsilon(t)\lvert^{2}y_{k}(t)dx+O(e^{-2\theta_{0}t}).
\end{eqnarray*}
From the result of Lemma 5.10, we have
\begin{eqnarray*}
\rvert I_{k}(t)-I_{k}(T_{n})\lvert\leq \frac{CA_{0}^{2}}{L}e^{-2\theta_{0}t}.
\end{eqnarray*}
Thus, by $\omega_{k}(T_{n})=\omega_{k}^{0}$ and $\varepsilon(T_{n})\equiv 0$, we obtain
\begin{eqnarray}
\rvert\int Q^{2}_{\omega_{k}(t)}dx-\int Q^{2}_{\omega_{k}^{0}}dx\lvert\leq C\|\varepsilon(t)\|^{2}_{L^{2}}+C(\frac{A_{0}^{2}}{L}+1)e^{-2\theta_{0}t}.
\end{eqnarray}
Recall that $\frac{d}{d\omega}\int Q^{2}_{\omega}dx\big\rvert_{\omega=\omega_{k}^{0}}>0$, then we assume $\omega_{k}(t)$ is close to $\omega_{k}^{0}$. Thus
\begin{align*}
(\omega_{k}(t)-\omega^{0}_{k})(\frac{d}{d\omega}\int Q^{2}_{\omega}dx\rvert_{\omega=\omega^{0}_{k}})
=&\int Q^{2}_{\omega_{k}(t)}dx-\int Q^{2}_{\omega_{k}^{0}}dx\notag\\
-&\beta(\omega_{k}(t)-\omega^{0}_{k})(\omega_{k}(t)-\omega^{0}_{k})^{2}
\end{align*}
with $\beta(\epsilon)\rightarrow 0$ as $\epsilon\rightarrow 0$, which implies that for some constant $C=C(\omega_{k}^{0})$,
\begin{eqnarray}
\rvert\omega_{k}(t)-\omega^{0}_{k}\lvert\leq C \rvert\int Q^{2}_{\omega_{k}(t)}dx-\int Q^{2}_{\omega_{k}^{0}}dx\lvert.
\end{eqnarray}
Therefore by (5.38) and (5.39), we have
\begin{eqnarray*}
\rvert\omega_{k}(t)-\omega^{0}_{k}\lvert\leq C\|\varepsilon(t)\|^{2}_{L^{2}}+C(\frac{A_{0}^{2}}{L}+1)e^{-2\theta_{0}t}.
\end{eqnarray*}
\par  This proves Claim 5.11.
\par  The following lemma gives a coercivity property of $H_{\omega_{k}^{0}}$.
\par \noindent\textbf{Lemma 5.12.} Let $d\geq 2,\;\;1+\frac{4}{d} < p < \frac{d+2}{(d-2)^{+}}$ and $\omega_{k}^{0}\in (\omega_{q},\;\omega_{p})$. Then there exists $\lambda > 0$ such that for any real-valued $v\in H^{1}$ satisfying $Re(Q_{\omega_{k}^{0}},\;v)=Im(Q_{\omega_{k}^{0}},\;v)=0$ and $Re(\nabla Q_{\omega_{k}^{0}},\;v)=0$, one has that
\begin{eqnarray}
(H_{\omega_{k}^{0}}v,v)\geq \lambda\|v\|_{H^{1}}^{2},
\end{eqnarray}
where $H_{\omega_{k}^{0}}$ is defined as (4.30) with $\omega=\omega_{k}^{0}$.
\par  \noindent\textbf{Proof.}  By (4.37) and Theorem 4.6, we have that $D''(\omega^{0}_{k}) > 0$. From Theorem 3.3 and Corollary 3.31 in \cite{GSS1987}, we get this result.

\par \noindent\textbf{Lemma 5.13.} Let $d\geq 2$ and $1+\frac{4}{d} < p < \frac{d+2}{(d-2)^{+}}$. For $\omega_{k}^{0}\in (\omega_{q},\;\omega_{p})$, $\omega_{k}(t)$ close to $\omega_{k}^{0}$ and $\Gamma_{\omega_{k}^{0}}(z)=E(z)+\omega_{k}^{0}M(z)$, we have
\begin{eqnarray}
\rvert\Gamma_{\omega_{k}^{0}}(Q_{\omega_{k}(t)})-\Gamma_{\omega_{k}^{0}}(Q_{\omega_{k}^{0}})\lvert\leq C\rvert\omega_{k}(t)-\omega_{k}^{0}\lvert^{2}.
\end{eqnarray}

\par \noindent\textbf{Proof.} By (2.2) and (2.3), we have
\begin{align}
\Gamma_{\omega_{k}^{0}}(Q_{\omega_{k}(t)})=E(Q_{\omega_{k}(t)})+\omega_{k}^{0}\int \rvert Q_{\omega_{k}(t)}\lvert^{2}dx.
\end{align}
By Taylar expansion of $\Gamma_{\omega_{k}^{0}}(Q_{\omega_{k}(t)})$, (5.42), Theorem 4.3 and Theorem 4.6,
\begin{align}
\Gamma_{\omega_{k}^{0}}(Q_{\omega_{k}(t)})
={} &\Gamma_{\omega_{k}^{0}}(Q_{\omega_{k}^{0}})-(\omega_{k}(t)-\omega_{k}^{0})^{2}\frac{d}{d\omega}\int Q^{2}_{\omega}dx\big\rvert_{\omega=\omega_{k}^{0}}\notag\\
+&\rvert\omega_{k}(t)-\omega_{k}^{0}\lvert^{2}\beta(\rvert\omega_{k}(t)-\omega_{k}^{0}\lvert).
\end{align}
By (5.43), Theorem 4.6 and $\omega_{k}(t)$ close to $\omega^{0}_{k}$, there exists $C=C(\omega^{0}_{k}) > 0$ such that
\begin{eqnarray*}
\rvert\Gamma_{\omega_{k}^{0}}(Q_{\omega_{k}(t)})-\Gamma_{\omega_{k}^{0}}(Q_{\omega_{k}^{0}})\lvert\leq C\rvert\omega_{k}(t)-\omega_{k}^{0}\lvert^{2}.
\end{eqnarray*}
\par  This proves Lemma 5.13.
\par  Now we set
\begin{equation}
J(t)=\sum_{k=1}^{K}[(\omega^{0}_{k}+\frac{1}{4}\rvert v_{k}\lvert^{2})I_{k}(t)-v_{k}M_{k}(t)]
\end{equation}
and
\begin{equation}
G(t)=E(\varphi_{n}(t))+J(t).
\end{equation}
By (5.33), (5.34), (5.35), Lemma 2.4 and 2.5 as well as Lemma 5.12 and 5.13, we directly deduce the following lemmas.
\par \noindent\textbf{Lemma 5.14.}  (Lemma 6(i) in \cite{MM2006})
For all $t\in[t^{*},T_{n}]$, we have
\begin{align}
G(t)&=\sum_{k=1}^{K}[E(Q_{\omega_{k}^{0}})+\omega_{k}^{0}\int Q^{2}_{\omega_{k}^{0}}dx]+P(\varepsilon(t),\varepsilon(t))+\sum_{k=1}^{K}O(\rvert\omega_{k}(t)-\omega_{k}^{0}\lvert^{2})\notag\\
&+\|\varepsilon(t)\|_{H^{1}}^{2}\beta (\|\varepsilon(t)\|_{H^{1}})+O(e^{-2\theta_{0}t})
\end{align}
with $\beta(\epsilon)\rightarrow 0$ as $\epsilon\rightarrow 0$, where
\begin{align}
P(\varepsilon(t),\varepsilon(t))={} &\int\rvert\nabla\varepsilon(t)\lvert^{2}dx-\sum_{k=1}^{K}[\int\rvert\widetilde{R}_{k}(t)\lvert^{\frac{4}{d}}\rvert\varepsilon(t)\lvert^{2}+\frac{4}{d}\rvert\widetilde{R}_{k}(t)\lvert^{\frac{4}{d}-2}[Re( \overline{\widetilde{R}}_{k}(t)\varepsilon(t))]^{2}dx]\notag\\
+&\sum_{k=1}^{K}[(\omega_{k}(t)+\frac{1}{4}\rvert v_{k}\lvert^{2})\int \rvert\varepsilon(t)\lvert^{2} y_{k}(t)dx-v_{k}\cdot Im\int \nabla \varepsilon(t)\cdot\overline{\varepsilon}(t)y_{k}(t)dx]\notag\\
+&\sum_{k=1}^{K}[\int\rvert\widetilde{R}_{k}(t)\lvert^{p-1}\rvert\varepsilon(t)\lvert^{2}+(p-1)\rvert\widetilde{R}_{k}(t)\lvert^{p-3}[Re( \overline{\widetilde{R}}_{k}(t)\varepsilon(t))]^{2}dx].
\end{align}

\par \noindent\textbf{Lemma 5.15.} (Lemma 4.1 in \cite{MMT2006})There exists $\lambda>0$ such that
\begin{eqnarray}
P(\varepsilon(t),\varepsilon(t))\geq \lambda\|\varepsilon(t)\|_{H^{1}}^{2},~~~~t\in[t^{*},T_{n}].
\end{eqnarray}
\par  Combining with Lemma 5.10, Lemma 5.11, Lemma 5.14 and Lemma 5.15, we can deduce the following lemma according to Martel and Merle's way \cite{MM2006}.
\par \noindent\textbf{Lemma 5.16.} (Lemma 5 in \cite{MMT2006})For any $t\in[t^{*},T_{n}]$,
\begin{eqnarray}
\|\varepsilon(t)\|_{H^{1}}^{2}+\rvert\omega_{k}(t)-\omega_{k}^{0}\lvert+\rvert x_{k}(t)\lvert^{2}+\rvert\gamma_{k}(t)-\gamma_{k}^{0}\lvert^{2}\leq C(\frac{A_{0}^{2}}{L}+1)e^{-2\theta_{0}t}.
\end{eqnarray}
\par \noindent\textbf{Lemma 5.17.}  For any $t\in [t^{*},T_{n}]$, there exists $C > 0$ such that
\begin{equation}
\|R(t)-\tilde{R}(t)\|_{H^{1}}\leq C\sum_{k=1}^{K}(\rvert\omega_{k}(t)-\omega_{k}^{0}\lvert+\rvert x_{k}(t)\lvert+\rvert\gamma_{k}(t)-\gamma_{k}^{0}\lvert).
\end{equation}
\noindent\textbf{Proof.} By (5.5) and (5.23),
\begin{align}
\widetilde{R}_{k}(t)-R_{k}(t)=&Q_{\omega_{k}(t)}(x-\widetilde{x}_{k}(t))e^{i(\frac{1}{2}v_{k}x+\delta_{k}(t))}-Q_{\omega^{0}_{k}}(x-\widetilde{x}^{0}_{k})e^{i(\frac{1}{2}v_{k}x+\delta^{0}_{k})}\notag\\
=&\nabla Q_{\omega_{k}(t)}(\eta_{k})e^{i(\frac{1}{2}v_{k}x+\delta_{k}(t))}x_{k}(t)
+iQ_{\omega_{k}(t)}(x-\widetilde{x}_{k}(t))e^{i\zeta_{k}}(\gamma_{k}(t)-\gamma^{0}_{k})\notag\\
+&
\frac{\partial Q_{\omega_{k}(t)}}{\partial \omega}(x-\widetilde{x}_{k}(t))\rvert_{\omega=\xi_{k}}e^{i(\frac{1}{2}v_{k}x+\delta_{k}(t))}(\omega_{k}(t)-\omega^{0}_{k}),
\end{align}
where $\eta_{k}$ is between $x-\widetilde{x}^{0}_{k}$ and $x-\widetilde{x}_{k}(t)$, $\zeta_{k}$ is between $\gamma^{0}_{k}$ and $\gamma_{k}(t)$, $\xi_{k}$ is between $\omega^{0}_{k}$ and $\omega_{k}(t)$. Thus
\begin{align*}
\|R(t)-\tilde{R}(t)\|_{H^{1}}\leq C\sum_{k=1}^{K}(\rvert\omega_{k}(t)-\omega_{k}^{0}\lvert+\rvert x_{k}(t)\lvert+\rvert\gamma_{k}(t)-\gamma_{k}^{0}\lvert).
\end{align*}

\par  This completes the proof of Lemma 5.17.

\noindent\textbf{Proof of Proposition 5.7.} From Lemma 5.16, we get for all $t\in[t^{*},T_{n}]$
\begin{align}
\|R(t)-\tilde{R}(t)\|_{H^{1}}^{2}\leq& C\sum_{k=1}^{K}(\rvert\omega_{k}(t)-\omega_{k}^{0}\lvert^{2}+\rvert\gamma_{k}(t)-\gamma_{k}^{0}\lvert^{2}+\rvert x_{k}(t)\rvert^{2})\notag\\
\leq& C(\frac{A_{0}^{2}}{L}+1)e^{-2\theta_{0}t},
\end{align}
By Lemma 5.17 and (5.52), we have
\begin{eqnarray*}
\|\varphi_{n}(t)-R(t)\|_{H^{1}}^{2}\leq 2\|\varepsilon(t)\|_{H^{1}}^{2}+2\|\tilde{R}(t)-R(t)\|_{H^{1}}^{2}\leq C(\frac{A_{0}^{2}}{L}+1)e^{-2\theta_{0}t},
\end{eqnarray*}
where $C>0$ does not depend on $A_{0}$. Now we choose $A_{0}^{2}>8C,\;L=A_{0}^{2}$, and $T_{0}$ large enough. It follows that
\begin{eqnarray*}
\|\varphi_{n}(t)-R(t)\|_{H^{1}}^{2}\leq 2Ce^{-2\theta_{0}t}\leq \frac{A_{0}^{2}}{4}e^{-2\theta_{0}t}.
\end{eqnarray*}
Therefore, the conclusion is that for any $t\in[t^{*},T_{n}]$,
\begin{align*}
\|\varphi_{n}(t)-R(t)\|_{H^{1}}\leq \frac{A_{0}}{2}e^{-\theta_{0}t}.
\end{align*}
\par  This completes the proof of Proposition 5.7.

\par \noindent\textbf{Acknowledgment.}
\par  This research is supported by the National Natural Science Foundation of China 11871138.

{\small


\begin{thebibliography}{50}
\bibitem{BJS2016}Bartsch, T., Jeanjean, L., Soave, N. (2016). Normalized solutions for a system of coupled cubic Schr\"odinger equations on $\mathbb{R}^{3}$. J. Math. Pures Appl (9), 106: 583-614.
\bibitem{BMRV2021}Bartsch, T., Molle, R., Rizzi, M., Verzini, G., (2021). Normalized solutions of mass supercritical Schr\"odinger equations with potential.  Comm. Partial Differential Equations  46(9): 1729-1756.
\bibitem{BZZ2020} Bartsch, T., Zhong, X., Zou, W. M. 2020. Normalized solutions for a system of coupled Schr\"odinger system. Math. Ann. 380(3-4): 1713-1740.
\bibitem{BC1981}Berestycki, H., Cazenave, T. (1981). Instabilit$\acute{e}$ des $\acute{e}$tats stationnaies dans les $\acute{e}$quations de Schr\"odinger et de Klein-Gordon non lin$\acute{e}$aires. C. R. Acad. Sci. Paris S$\acute{e}$r. I Math. 293: 489-492.
\bibitem{BL1983}Berestycki, H., Lions, P. L. (1983). Nonlinear scalar field equations. Arch. Rational Math. Anal., 82: 313-375.
\bibitem{C2003}Cazenave, T. (2003). Semilinear Schr\"odinger equations. Courant Lecture Notes in Mathematics, 10. New York University, Courant Institute of Mathematical Sciences, New York; American Mathematical Society, Providence, RI.
\bibitem{CL1982}Cazenave, T., Lions, P. L. (1982).  Orbital stability of standing waves for some nonlinear Schr\"odinger equations. Comm.  Math. Phys. 85(4): 549-561.
\bibitem{CW1990}Cazenave, T., Weissler, F. B. (1990).  The Cauchy problem for the critical nonlinear Schr\"odinger equation in $H^{s}$, Nonlinear Anal. 14(10): 807-836.
\bibitem{CC2011}C$\hat{o}$te, R., Le Coz, S. (2011). High-speed excited multi-solitons in nonlinear Schr\"odinger equations. J. Math. Pures Appl. 96: 135-166.
\bibitem{CV2017}Cirant, M., Verzini, G. (2017). Bifurcation and segregation in quadratic two-populations mean field games systems.
ESAIM Control Optim. Calc. Var. 23, no. 3, 1145-1177.
\bibitem{CMM2011}C$\hat{o}$te, R., Martel, Y., Merle, F. (2011). Construction of multi-soliton solutions for the $L^{2}$-supercritical gkdv and NLS equations. Revista Matem$\acute{a}$tica Iberomamericana. 27(1): 273-302.
\bibitem{FH2021}Fukaya, N., Hayashi, M. (2021). Instability of algebraic standing waves for nonlinear Schr\"odinger equations with double power nonlinearities. Transactions of the American Mathhematical Society. 374: 1421-1447.
\bibitem{F2003}Fukuizumi, R. (2003). Stability and instability of standing waves for nonlinear Schr\"odinger equations. Tohoku Mathematical Publications. No. 25.
\bibitem{GNN1979}Gidas, B., Ni, W. M., Nirenberg, L. (1979). Symmetry and related properties via the maximal principle. Comm. Math. Phys. 68: 209-243.
\bibitem{GSS1987}Grillakis, M., Shatah, J., Strauss, W. A. (1987). Stability theory of solitary waves in the presence of symmetry, I. J. Funct. Anal. 74(1): 160-197.
\bibitem{GSS1990}Grillakis, M., Shatah, J., Strauss, W. A. (1990). Stability theory of solitary waves in the presence of symmetry,II. J. Funct. Anal. 94(2): 308-348.
\bibitem{K1959}Kato, T. (1959). Growth properties of solutions of the reduced wave equation with a variable coefficient. Comm, Pure Appl. Math. 12: 403-425.
\bibitem{K2011}Kawano, S. (2011). Uniqueness of positive solutions to semilinear elliptic equations with double power nonlinearities. Differential  Integral Equations. 24(24): 201-207.
\bibitem{K1989}Kwong, M. K. (1989). Uniqueness of positive solutions of $\triangle u-u+u^{p}=0$ in $\mathbb{R}^{n}$. Arch. Rational Mech. Anal. 105(3): 243-266.
\bibitem{CMR2016}Le Coz, S., Martel, Y., Rapha\"el, P. (2016) Minimal mass blow up solutions for a double power nonlinear Schr\"odinger equation. Revista Matem$\acute{a}$tica Iberomamericana. 32(3): 795-833.
\bibitem{LL2000}Lieb, E. H., Loss, M. (2000). Analysis, Graduate Studies in Mathematics Volume 14, American Mathematical Society.
\bibitem{MM2006}Martel, Y., Merle, F. (2006). Multi solitary waves for nonlinear Schr\"odinger equations, Ann. I. H. Poincar\'{e}-AN. 23: 849-864.
\bibitem{MMT2006}Martel, Y., Merle, F., Tsai, T. P. (2006). Stability in $H^1$ of the sum of $K$ solitary waves for some nonlinear Schr\"odinger equations. Duke Mathematical Journal. 133(3): 405-466.
\bibitem{M1993}McLeod, K., (1993). Uniqueness of positive radial solutions of $\Delta u+f(u)=0$ in $\mathbb{R}^{n}$. II, Trans. Amer. Math. Soc. 339 495-505.
\bibitem{M1990}Merle, F. (1990). Construction of solutions with exactly k blow-up points for the Schr\"odinger equation with critical nonlinearity. Comm. Math. Phys. 129 (2): 223-240.
\bibitem{NS2012}Nakanishi, K.,  Schlag, W. (2012). Global dynamics above the ground state energy for the cubic NLS equation in 3D. Calculus of Variations and Partial Differential Equations. 44(1-2): 1-45.
\bibitem{P1965}Pohozaev, S. I. (1965). Eingenfunctions of the equations of the $\Delta u+\lambda f(u)=0$. Sov. Math. Doklady. 165: 1408-1411.
\bibitem{RS1978}Reed, M., Simon, B. (1978). Methods of Modern Mathematical Physics. Vol. IV, Academic Press,  New York.
\bibitem{S2009} Schlag, W. (2009). Stable manifolds for an orbitally unstable NLS. Ann. of Math. (2) 169, no. 1: 139-227.
\bibitem{S2020}Soave, N. (2020). Normalized ground states for the NLS equation with combined nonlinearities. J. Differential Equations. 269: 6941-6987.
\bibitem{S1977}Strauss, W. A. (1977). Existence of solitary waves in higher dimensions. Comm. Math. Phys. 55(2): 149-162.
\bibitem{TVZ2007}Tao, T., Visan, M., Zhang, X. (2007).The nonlinear schrdinger equation with combined power-type nonlinearities. Comm. Partial Differential Equations, 32(8): 1281-1343.
\bibitem{T1987} Tsutsumi, Y. (1987). $L^{2}$-solutions for nonlinear Schr\"odinger equations and nonlinear group, Funkcial. Ekvac. 30: 115-125.
\bibitem{T2009}Tao, T. (2009). Why are solitons stable?. Bulletin of the American Mathematical Society. 46(1): 1-33.
\bibitem{W1983}Weinstein, M. I. (1983). Nonlinear Schr\"odinger equations and sharp interpolations estimates. Comm. Math. Phys. 87: 567-576.
\bibitem{W1985}Weinstein, M. I. (1985). Modulation stability of ground states of nonlinear Schr\"odinger equations. SIAM J. Math. Anal. 16: 472-491.
\bibitem{W1986}Weinstein, M. I. (1986). Lyapunov stability of ground states of nonlinear dispersive evolution equations. Comm. Pure  Appl. Math.  39: 51-68.
\bibitem{Z2000}Zhang, J. (2000). Stability of attractive Bose-Einstein condensates, Journal of Statistical Physics. (3-4) 101: 731-746.
\bibitem{Z2005}Zhang, J. (2005). Sharp threshold for blowup and global existence in nonlinear Schr\"odinger equations under a harmonic potential. Comm. Partial Differential Equations 30: 1429-1443.



\end{thebibliography}
\end{document}